\documentclass[english]{article}
\usepackage[utf8]{inputenc}
\usepackage[T1]{fontenc}      
\usepackage[english]{babel}
\usepackage{graphicx}
\usepackage{amsmath,amsthm}
\usepackage{amsfonts,amssymb}
\usepackage{lmodern}
\usepackage[a4paper]{geometry}
\usepackage{mathtools}
\usepackage{array}
\usepackage[toc,page]{appendix}
\usepackage{abstract}
\usepackage{caption}
\usepackage{subcaption}
\usepackage{multirow}
\usepackage{authblk} 
\newcommand*\samethanks[1][\value{footnote}]{\footnotemark[#1]}

\begin{document}
\title{\textbf{Exact simulation of the jump times of a class of Piecewise Deterministic Markov Processes}}

\author{\textbf{Vincent Lemaire}\thanks{Laboratoire de Probabilités et Modèles Aléatoires, UMR CNRS 7599, UPMC Paris 6, Case 188, 4 place Jussieu, F-75252 Paris Cedex 5, France} \footnote{vincent.lemaire@upmc.fr}\hspace{1cm}\textbf{Michèle Thieullen}\samethanks[1] \footnote{michele.thieullen@upmc.fr}  \hspace{1cm}\textbf{Nicolas Thomas}\samethanks[1] \footnote{nicolas.thomas@upmc.fr}}

\date{}

\maketitle

\begin{abstract}
\noindent
In this paper, we are interested in the exact simulation of a class of Piecewise Deterministic Markov Processes (PDMP). We show how to perform efficient thinning algorithms depending on the jump rate bound. For different types of jump rate bounds, we compare theoretically the mean number of generated (total) jump times and we compare numerically the simulation times. We use thinning algorithms on Hodgkin-Huxley models with Markovian ion channels dynamic to illustrate the results. 
\end{abstract}

\vspace{0.2cm}
\noindent \textbf{Keywords:} Piecewise Deterministic Markov Process; Exact simulation; Thinning; Neural model; Channel noise 
\vspace{0.2cm}\\
\textbf{Mathematics Subject Classification:} 60J25, 65C20, 92C20, 92C45 

\section{Introduction}
In many areas it is important to be able to simulate exactly and rapidly trajectories of a stochastic process. This is the case for Monte Carlo methods, statistical estimation, bootstrap. In this article, we are interested in the exact simulation (perfect sampling) of a class of Piecewise Deterministic Markov Processes (PDMP). These processes, introduced by M.H.A. Davis in \cite{davarticle}, are based on an increasing sequence of random times in which the processes have a jump and on a deterministic evolution between two successive random times. The law of a PDMP is thus determined by three parameters called  the characteristics of the PDMP: a family of vector fields, a jump rate (intensity function) and a transition measure. These parameters can be given but in some cases we may have only access to few of them. In this study we suppose that the flows of the PDMP are known, this means that we explicitly know the solution of each ordinary differential equations associated to each vector fields. This hypothesis is not a restriction. Indeed, many random phenomena are modelled by PDMPs with explicit flows. For example, we can quote the temporal evolution of the membrane potential and ionic channels in neuroscience (see \cite{wainfluid}), the evolution of a food contaminant in the human body in pharmacokinetic (see \cite{bouguet}), the growth of bacteria in biology (see \cite{doumic}), the data transmission in internet network (see \cite{malri}) or the evolution of a chemical network in chemistry (see \cite{alfonsi}).\\
The flows being known, we focus on the simulation of the jump times. The basic tool we use is the thinning method introduced by Lewis and Schedler in \cite{lew} to simulate Poisson processes and generalised by Ogata \cite{ogata} to simulate any point processes. This method became classic when the jump rate of the process admits a constant upper bound $\overline{\lambda}$. It consists in generating the jump times of a (homogeneous) Poisson process with intensity $\overline{\lambda}$ then selecting some of these times by a rejection argument. In this case, we obtain an algorithm that is easy to implement. However, it is intuitive that a constant bound $\overline{\lambda}$ will not provide a powerful algorithm when the jump rate of the process present significant variations.\\
Let us go back to PDMPs. The main result of this article is to develop thinning algorithms where the jump rate bound closely follows the evolution in time of the jump rate of the PDMP to simulate. In the sequel, such a bound is called \textit{optimal bound}.\\
We see at least three interests in this result. The first is that it applies with weak hypotheses on the jump rate: we are not assuming that the jump rate is bounded, it can depend on time, being monotone or not. The second is that it provides powerful thinning algorithms.
The drawback of this method is that more the jump rate bound is close to the jump rate more the method is time consuming. It is thus necessary to look for a satisfactory balance. We discuss this difficulty on a numerical example. Finally, the \textit{optimal bound} is constructed by following each vector fields of the family. This construction is thus natural in the context of switching processes such as PDMPs. For this reason we think that the algorithms studied in this article can be applied to a much larger family of processes such as Hawkes processes or jump (P)SDE with state-dependent intensity.\\
Let us now give some details on the content of the article. In this study we will consider three types of jump rate bounds:
\begin{itemize}
\item the \textit{global bound}, the coarsest, which is constant (in particular it is independent of the state of the PDMP and of time),
\item the \textit{local bound}, which depends on each post-jump values of the PDMP and which is constant between two successive jump times,
\item the \textit{optimal bound}, the finest, which depends on each post-jump values of the PDMP and also on the time evolution of the process between two successive jump times.
\end{itemize}    
For each of the three jump rate bounds we theoretically compare the mean number of rejected points, then, on a numerical example, the simulation times.\\
As an indicator of the efficiency of the thinning method, we choose the mean value of the ratio between the number of accepted jump times and the number of proposed jump times. We call it rate of acceptance. This indicator is between $0$ and $1$ and is easily understood, the closer it is to $1$ the less we reject points, thus the more the method is efficient. We explicitly express this rate of acceptance in terms of the transition measure of a discrete time Markov chain which carries more than the information of the PDMP only. Note that it is not the embedded Markov chain of the PDMP since it contains also all the rejected jump times. The rate of acceptance is also expressed as a function of the ratio between the jump rate of the PDMP and the jump rate bound. This expression allows to see that more the jump rate bound is close to the jump rate of the PDMP more the algorithm is efficient. 
Let us also note that the rate of acceptance is different from the definition of efficiency given in \cite{lew} or \cite{dev} chap. 6 where it is defined as the ratio between the mean number of accepted jump times and the mean number of proposed jump times. However, in the case of Poisson processes, the rate of acceptance coincides with the definition of efficiency given in \cite{lew} or \cite{dev} chap. 6.\\
When the flows are not known explicitly, we can simulate PDMPs with the algorithms given in \cite{rid} or in \cite{veltz} which essentially consist in repeatedly solving ODEs. \\
We chose to illustrate the results on the efficiency of the \textit{optimal bound} on two stochastic Hodgkin-Huxley (HH) models, the \textit{subunit model} and the \textit{channel model} \cite{wainfluid}. The \textit{subunit model} is the stochastic model whose deterministic limit (when the number of channels goes to infinity) is the classical four-dimensional deterministic HH model \cite{hodgkin} whereas the deterministic limit of the \textit{channel model} is a fourteen-dimensional model \cite{wainfluid}. The difference between the two stochastic models comes from the modelling of the conductance of the membrane. In the \textit{subunit model}, the conductance is given by the proportion of open gates whereas, in the \textit{channel model}, it is given by the proportion of open channels. Note that the two deterministic models are equivalent when the initial conditions satisfy a binomial relation \cite{wainfluid}. The \textit{channel model} is known to be more bio-realistic than the \textit{subunit model} but also more complex.\\
Thinning algorithms are better suited than algorithms in \cite{rid} or\cite{veltz} to simulate variables of biological interest such as mean spiking times or mean dwell times using classical Monte Carlo. Indeed, in addition to be faster, thinning algorithms do not introduce bias since the simulation is exact.\\
We present both stochastic HH models as PDMPs. The jump rates of these models present high variations especially when the membrane is in a depolarization phase, thus, it allows to check whether the \textit{optimal bound} speed up simulation compared to the \textit{global bound} and the \textit{local bound}. Moreover the jump rates (which come from the modelling \cite{hodgkin}) are complex functions, thus numerical inversion of the distribution function can be time consuming. To the best of our knowledge, when flows are known, no studies of the error have been carried out when we numerically inverse the distribution function.\\  
Several algorithms can be found to simulate stochastic HH models. These algorithms can be classified into three categories: channel-state tracking algorithms (CST) \cite{defelice}, \cite{rubin}, \cite{ander} channel-number-tracking algorithms (CNT) \cite{skaugen}, \cite{chowwhite}, \cite{ander} and approximate algorithms \cite{soudry}, \cite{gold}, \cite{fox}. The main difference between the CNT and CST algorithms comes from the simulation of the inter-jump times and we emphasise that thinning algorithms can be used in both types of algorithm.\\
In the literature, the term exact algorithm (simulation) is employed to denote CST and CNT algorithms (see \cite{fox}, \cite{ander}) even if we use some Euler integration to numerically solve the inversion of the distribution function problem. In this article, the meaning of "exact" signifies that the simulated trajectories are realizations of these
processes.\\
We show how to determine jump rate bounds in such stochastic HH models. We also present a new way to simulate the transition measure of the \textit{channel model} by using the transition measure of the \textit{subunit model}. This approach is numerically more efficient  than the classical method used for example in \cite{mino} p.587. For the stochastic HH models studied here, we show that the simulation time is reduced by 2 in going from the \textit{global bound} to the \textit{local bound} and that it is also reduced by 2 in going from the \textit{local bound} to the \textit{optimal bound}.\\
The paper is organised as follows. In section 2, we give the definition of PDMPs and set the
notations used in other sections. In section 3, we introduce the jump rate bounds and we
present the thinning procedure. In section 4, we give the theoretical results concerning the comparison of the jump rate bounds and the rate of acceptance. In section 5, we introduce the Hodgkin-Huxley models. In section 6 we numerically illustrates the results.    \\   

\section{Formal definition of PDMPs}
A PDMP is a stochastic process in which the source of randomness comes from random jump times and post-jump locations. In this paper, we consider that such a process takes the following general form
\[x_{t}=(\theta_{t},V_{t})\hspace{0,5cm}\forall{t}\geq0
\]
where:

\begin{itemize}

	\item $\theta :\mathbb{R}^{+}\to{K}$ is a jump process that characterizes the mode of the system, $m\geq1$.
	
	\item $V :\mathbb{R}^{+}\to{D}\subseteq{\mathbb{R}^{d}}$ is a stochastic process which evolves deterministically between two successive jumps of the process $\theta$, $d\geq1$.
	
	\item $K$ is a finite or countable space.
	
\end{itemize}

\noindent
Let us denote $E=K\times{D}$ so that $(x_{t})_{t\geq0}$ is an $E$-valued process. We note $(T_{n})_{n\geq0}$ the sequence of jump times of the PDMP and $(N_{t})_{t\geq0}$ the counting process $N_{t}=\sum_{n\geq1}\textbf{1}_{T_{n}\leq{t}}$. We make the following assumption
\vspace{0.2cm}\\
\textbf{\underline{Assumption 2.1}} : for every starting point $x\in{E}$, $\mathbb{E}_{x}[N_{t}]<\infty$ for all $t\geq0$.
\vspace{0.2cm}\\
This assumption implies in particular that $T_{n}\rightarrow\infty$ almost surely.
Such a PDMP is uniquely defined by its three characteristics $(\phi,\lambda,Q)$ where

\begin{itemize}

	\item \textit{The deterministic flow} $\phi:\mathbb{R}_{+}\times{E}\to{D}$ is supposed continuous and induced by a vector field $F : E\rightarrow{D}$.
For $t\in{[T_{n},T_{n+1}[}$, $V$ takes the following form
$
V_{t}=\phi(t-T_{n},x_{T_{n}})
$ and the trajectory of the process $(x_{t})_{t\geq0}$ is then given by
	\[x_{t}=\sum_{n\geq0}\Big{(}\theta_{T_{n}},\phi(t-T_{n},x_{T_{n}})\Big{)}\textbf{1}_{T_{n}\leq{t}<T_{n+1}}.\]
	For notational convenience, we define a vector field $G : E\rightarrow{E}$ such that for $x\in{E}$ $G(x)=\left(\begin{array}{c}0\\F(x)
	\end{array}\right)$. Then we can represent the PDMP as follows 
	\[x_{t}=\sum_{n\geq0}\psi(t-T_{n},x_{T_{n}})\textbf{1}_{T_{n}\leq{t}<T_{n+1}}.\]
	Where $\psi$ is the flow induced by $G$ such that for $t\in{[T_{n},+\infty[}$ \[\psi(t-T_{n},x_{T_{n}})=\Big{(}\theta_{T_{n}},\phi(t-T_{n},x_{T_{n}})\Big{)}.\]
	
	\item \textit{The jump rate} $\lambda:E\to\mathbb{R}_{+}$ is a non-negative measurable function that characterizes the frequency of jumps and such that for each $x\in{E}$ there exists $\sigma(x)$ such that the function $s\rightarrow\lambda\Big{(}\psi(s,x)\Big{)}$ is integrable on $[0,\sigma(x)[$. Then the trajectory of the (stochastic) jump rate is given by 
	  	\[\lambda(x_{t})=\sum_{n\geq0}\lambda(\psi(t-T_{n},x_{T_{n}}))\textbf{1}_{T_{n}\leq{t}<T_{n+1}}.\]
	  	The formula is to be understood as follows: the intensity of the inter-jump time $S_{n+1}$ is $\lambda(\psi(t-T_{n},x_{T_{n}}))$ for $t\geq{T_{n}}$. 
	\item \textit{The transition measure} $Q:E\times{\mathcal{B}(E)}\to[0,1]$ governs the post-jump location of the process. It verifies
	\[Q\Big{(}x,\{x\}\Big{)}=0\hspace{1cm}\forall{x\in{E}}.
	\]
	
	\end{itemize}

\noindent	
The iterative construction of these processes in \cite{davarticle} and \cite{dav} provides a way to simulate their trajectories although the problem of the exact simulation of the inter-jump times is not obvious. In \cite{dav}, M.H.A Davis shows that there exists a filtered probability space $(\Omega,\mathcal{F}, \mathcal{F}_{t}, \mathbb{P}_{x})$ such that the process $(x_{t})_{t\geq0}$ is a Markov process. He also shows that $(x_{T_{k}})_{k\geq0}$ is a Markov chain with kernel $Z$ such that for $x_{0}\in{E}$
\begin{align*}
Z(x_{0},A)&=\mathbb{P}\Big{(}x_{T_{1}}\in{A}|x_{T_{0}}=x_{0}\Big{)}\\
&=\int_{0}^{\infty}Q\Big{(}\psi(t,x_{0})),A\Big{)}\lambda\Big{(}\psi(t,x_{0})\Big{)}
e^{-\int_{0}^{t}\lambda\Big{(}\psi(s,x_{0})\Big{)}ds}dt.
\end{align*}

\noindent
The randomness of the PDMP is contained in the associated jump process $(\eta_{t})$ defined by
\[\eta_{t}=x_{T_{n}}\hspace{1cm}T_{n}\leq{t}<T_{n+1}.\]
Because $T_{n}=\inf\{t>T_{n-1} : \theta_{t-}\neq{\theta_{t}}\}$, the knowledge of $(\eta_{t})_{t\geq0}$ implies the knowledge of $(T_{n})_{n\geq0}$.

\section{Simulation of PDMPs by thinning}

In this section, we first present the three different jump rate bounds. Secondly, we describe the thinning method to simulate inter-jump times of PDMPs. 

\subsection{Jump rate bounds} 

In this section we introduce the different jump rate bounds, namely, the \textit{optimal bound}, the \textit{local bound} and the \textit{global bound}. The \textit{optimal bound} is particularly efficient in term of reject because it is as close as we want to the jump rate. Let us first introduce a general bound of $\lambda$, namely $\tilde{\lambda}$, defined by
\[\tilde{\lambda}(x_{t},t)=\sum_{n\geq0}\tilde{\lambda}(\psi(t-T_{n},x_{T_{n}}),t)\textbf{1}_{T_{n}\leq{t}<T_{n+1}}\]
where $(T_{n})_{n\geq0}$ denotes the jump times of the PDMP. For $n\geq0$ and conditionally on $(T_{n},x_{T_{n}})$ the function $t\longmapsto\tilde{\lambda}(\psi(t-T_{n},x_{T_{n}}),t)$ is defined on $[T_{n},+\infty[$ and verifies
\[\lambda(\psi(t-T_{n},x_{T_{n}}))\leq{}\tilde{\lambda}(\psi(t-T_{n},x_{T_{n}}),t)\hspace{1cm}\forall{t}\geq{T_{n}}.\]
We shall precise that the function of time $\tilde{\lambda}(\psi(t-T_{n},x_{T_{n}}),t)$ is used to simulate $S_{n+1}$.
For $n\geq0$ and conditionally on $(T_{n},x_{T_{n}})$ we define the function $\tilde{\Lambda}_{n}:\mathbb{R}_{+}\rightarrow\mathbb{R}_{+}$ by
\[\tilde{\Lambda}_{n}(t)=\int_{T_{n}}^{t}\tilde{\lambda}(\psi(s-T_{n},x_{T_{n}}),s)ds\]
We denote by $\left(\tilde{\Lambda}_{n}\right)^{-1}$  the inverse of $\tilde{\Lambda}_{n}$.

\subsubsection{The global bound}

We define the \textit{global bound} by 
\[\tilde{\lambda}^{\text{glo}}(x_{t},t)=\sum_{n\geq0}\sup_{x\in{E}}\lambda(x)\textbf{1}_{T_{n}\leq{t}<T_{n+1}}\]
By definition, this bound is constant and does not depend on the state of the PDMP nor on time, we will denote it by $\tilde{\lambda}^{\text{glo}}$. This bound is probably the most used and has the advantage to lead to an easy implementation. Indeed, to simulate the jump times of the PDMP we simulate a homogeneous Poisson process with jump rate $\tilde{\lambda}^{\text{glo}}$ disregarding the state of the PDMP.  For $n\geq0$ and conditionally on $(T_{n},x_{T_{n}})$, the integrated jump rate bound is given by $\tilde{\Lambda}^{\text{glo}}_{n}(t)=\tilde{\lambda}^{\text{glo}}(t-T_{n})$ for $t\geq{T_{n}}$ and the inverse is given by $(\tilde{\Lambda}_{n}^{\text{glo}})^{-1}(s)=
s/\tilde{\lambda}^{\text{glo}}+T_{n}$ for $s\geq0$. 

\subsubsection{The local bound}

We define the \textit{local bound} by 
\[\tilde{\lambda}^{\text{loc}}(x_{t},t)\equiv\sum_{n\geq0}\left(\sup_{t\geq{T_{n}}}\lambda\left(\psi(t-T_{n},x_{T_{n}})\right)\right)\textbf{1}_{T_{n}\leq{t}<T_{n+1}}\]
By definition, this bound is constant between two successive jump times and has the advantage of being adapted to the state of the PDMP right-after a jump. To each jump time of the PDMP corresponds a homogeneous Poisson process whose intensity depends on the state of the PDMP at the jump time. For $n\geq0$ and conditionally on $(T_{n},x_{T_{n}})$, the integrated jump rate bound is $\tilde{\Lambda}^{\text{loc}}_{n}(t)=\Big{(}\sup_{t\geq{T_{n}}}\lambda(\psi(t-T_{n},x_{T_{n}}))\Big{)}(t-T_{n})$ for $t\geq{T_{n}}$ and the inverse is given by $(\tilde{\Lambda}_{n}^{\text{loc}})^{-1}(s)=\Big{(}
s/\sup_{t\geq{T_{n}}}\lambda(\psi(t-T_{n},x_{T_{n}}))\Big{)}+T_{n}$ for $s\geq0$. 

\subsubsection{The optimal bound}

Let $P$ be a finite or a countable space. For $n\geq0$, let us denote by $(\mathcal{P}^{T_{n}}_{k})_{k\in{P}}$ a partition of $[T_{n},+\infty[$. Thus, for $k\in{P}$, $\mathcal{P}^{T_{n}}_{k}$ is an interval of $[T_{n},+\infty[$.
We define the \textit{optimal bound} by
\[\tilde{\lambda}^{\text{opt}}(x_{t},t)\equiv\sum_{n\geq0}\left(\sum_{k\in{P}}\sup_{t\in{\mathcal{P}^{T_{n}}_{k}}}
\lambda(\psi(t-T_{n},x_{T_{n}}))\textbf{1}_{\mathcal{P}^{T_{n}}_{k}}(t)\right)\textbf{1}_{T_{n}\leq{t}<T_{n+1}}\]
By definition, this bound is piecewise constant between two successive jump times, thus it is adapted to the state of the PDMP right-after a jump but also to the evolution in time of the jump rate. To each jump time of the PDMP corresponds a non-homogeneous Poisson process whose intensity depends on the state of the PDMP at the jump time and on the flow starting from this state.
For $n\geq0$ and conditionally on $(T_{n},x_{T_{n}})$, the integrated jump rate bound is, for $t\geq{T_{n}}$, 
\[
\tilde{\Lambda}^{\text{opt}}_{n}(t)=\sum_{k\in{P}}\sup_{t\in{\mathcal{P}^{T_{n}}_{k}}}\lambda\Big{(}\psi(t-T_{n},x_{T_{n}})\Big{)}\Big{|}\mathcal{P}^{T_{n}}_{k}\cap[T_{n},t]\Big{|}
\]
where $\Big{|}
\mathcal{P}^{T_{n}}_{k}\cap[T_{n},t]\Big{|}$ represents the length of the interval $\mathcal{P}^{T_{n}}_{k}\cap[T_{n},t]$.
As $P$ is at most countable, let us denote by $p_{i}$ for $i=0,\ldots,\text{Card}(P)$ its elements. The inverse, $\left(\tilde{\Lambda}^{\text{opt}}_{n}\right)^{-1}$, is given by
\[
\left(\tilde{\Lambda}^{\text{opt}}_{n}\right)^{-1}(s)=\sum_{i=0}^{\text{card}(P)}\Big{(}\frac{s-\sum_{
k=0}^{i-1}\sup_{t\in{\mathcal{P}^{T_{n}}_{p_{k}}}}
\lambda\Big{(}\psi(t-T_{n},x_{T_{n}})\Big{)}\Big{|}\mathcal{P}^{T_{n}}_{p_{k}}\Big{|}}{\sup_{t\in{\mathcal{P}^{T_{n}}_{p_{i}}}}
\lambda\Big{(}\psi(t-T_{n},x_{T_{n}})\Big{)}}+T_{n}+\sum_{l=0}^{i-1}\Big{|}\mathcal{P}^{T_{n}}_{p_{l}}\Big{|}\Big{)}
\textbf{1}_{[\kappa_{p_{i-1}},\kappa_{p_{i}}[}(s)
\]
where $\kappa_{p_{i}}=\sum_{
k=0}^{i}\sup_{t\in{\mathcal{P}^{T_{n}}_{k}}}
\lambda\Big{(}\psi(t-T_{n},x_{T_{n}})\Big{)}\Big{|}\mathcal{P}^{T_{n}}_{k}\Big{|}$ and, by convention $\sum_{l=0}^{-1}\Big{|}\mathcal{P}^{T_{n}}_{p_{l}}\Big{|}=0$ and $\kappa_{p_{-1}}=0$.
\vspace{0.3cm}\\
\underline{\textbf{Remark 1}}: Let $P=\mathbb{N}$ and let $\epsilon>0$. The \textit{optimal bound} with the partition $(\mathcal{P}_{k}^{T_{n},\epsilon})_{k\in{\mathbb{N}}}$ where, for $k\in{\mathbb{N}}$, $\mathcal{P}_{k}^{T_{n},\epsilon}=[T_{n}+k\epsilon,T_{n}+(k+1)\epsilon[$ can be use even if the jump rate $\lambda$ is not bounded. A bounded hypothesis on $\lambda$ is required to use the \textit{global bound} and it can be weakened for the \textit{local bound} (for example if the flows are bounded and $\lambda$ continuous).  
\vspace{0.2cm}\\
\underline{\textbf{Remark 2}}: For the three jump rate bounds, the simulation is exact. In particular, for all finite or countable $P$, that is, for any partitions of $[T_{n},+\infty[$, the simulation remains exact.
\vspace{0.2cm}\\
\underline{\textbf{Remark 3}}: The choice of the bound depends on the PDMP we want to simulate. If the jump rate does not vary very much in time, the \textit{local bound} or the \textit{global constant} bound can be adopted but if the jump rate presents high variations in a small time interval, the \textit{optimal bound} is preferable. 

\subsection{Thinning}

In this section, we detail the thinning procedure to simulate a PDMP.
Details on the thinning of Poisson processes may be find in \cite{lew} or \cite{dev}. We present the procedure with the general bound $\tilde{\lambda}$ of section 3.1. In practice, one has to apply the procedure with one of the three bounds (\textit{optimal}, \textit{local} or \textit{global}). We simulate a sample path of the PDMP $(x_{t})_{t\geq0}$ with values in $E$, starting from a fixed initial point $x_{T_{0}}=x_{0}$ at time $T_{0}$ as follows.\\
Let $(\tilde{T}^{0}_{k})_{k\geq0}$ be the Poisson process defined on $[0,+\infty[$ with jump rate $\tilde{\lambda}(\psi(t,x_{0}),t)$ for $t\geq0$, and
\[
\tau_{1}=\inf\{k>0 : U_{k}\tilde{\lambda}(\psi(\tilde{T}^{0}_{k},x_{0}),\tilde{T}^{0}_{k})\leq{\lambda}(\psi(\tilde{T}^{0}_{k},x_{0}))\}
\]
where $(U_{n})_{n\geq1}$ is a sequence of independent random variables with uniform distribution on $[0,1]$, independent of $(\tilde{T}^{0}_{k})_{k\geq0}$ and $x_{0}$.
The first jump time $T_{1}=S_{1}$ of the PDMP is the first jump of a non-homogeneous Poisson process defined on $[0,+\infty[$ with jump rate $\lambda(\psi(t,x_{0}))$ for $t\geq0$. We simulate $T_{1}$ by thinning the process $(\tilde{T}^{0}_{k})_{k\geq0}$, then, $T_{1}=\tilde{T}^{0}_{\tau_{1}}$. On $[0,T_{1}[$ the PDMP evolves as follows
\[x_{t}=\psi(t,x_{0})\]
and the random variable $x_{T_{1}}$ has distribution $Q\Big{(}(\psi(T_{1},x_{0})),.\Big{)}$. Note that conditionally on $(T_{0},x_{T_{0}},T_{1})$ the process $(\tilde{T}^{0}_{k})_{k\geq1}$ is a Poisson process on $[T_{0},T_{1}[$ with jump rate $\tilde{\lambda}(\psi(t,x_{0}),t)-\lambda(\psi(t,x_{0}))$, see \cite{dev} chap.6.
\vspace{0.2cm}\\
Suppose we have simulated $T_{i}$, then, conditionally on $(T_{i},x_{T_{i}})$, the PDMP $(x_{t})$ restarts from $x_{T_{i}}$ at time $T_{i}$ independently from the past. Let $(T_{i}+\tilde{T}^{i}_{k})_{k\geq0}$ be the Poisson process on $[T_{i},+\infty[$ with jump rate $\tilde{\lambda}(\psi(t-T_{i},x_{T_{i}}),t)$ for $t\geq{T_{i}}$ and 
\[
\tau_{i+1}=\inf\{k>0 : U_{k}\tilde{\lambda}(\psi(\tilde{T}^{i}_{k},x_{T_{i}}),T_{i}+\tilde{T}^{i}_{k})\leq{\lambda}(\psi(\tilde{T}^{i}_{k},x_{T_{i}}))\}
\]
where $(U_{n})_{n\geq1}$ is a sequence of independent uniform random variables, independent of $(\tilde{T}^{i}_{k})_{k\geq0}$ and $x_{T_{i}}$. $T_{i+1}$ is the first jump of a non-homogeneous Poisson process defined on $[T_{i},+\infty[$ with jump rate $\lambda(\psi(t-T_{i},x_{T_{i}}))$ for $t\geq{T_{i}}$. We simulate $T_{i+1}$ by thinning the process $(T_{i}+\tilde{T}^{i}_{k})_{k\geq0}$, then, $T_{i+1}=T_{i}+\tilde{T}^{i}_{\tau_{i+1}}$. On $[T_{i},T_{i+1}[$ the process evolves as follows
\[
x_{t}=\psi(t-T_{i},x_{T_{i}})
\]
and the random variable $x_{T_{i+1}}$ has distribution $Q\Big{(}(\psi(S_{i+1},x_{T_{i}})),.\Big{)}$.
Note that, conditionally on $(T_{i},x_{T_{i}},T_{i+1})$, the process $(T_{i}+\tilde{T}^{i}_{k})_{k\geq{1}}$ is a Poisson process on $[T_{i},T_{i+1}[$ with jump rate $\tilde{\lambda}(\psi(t-T_{i},x_{T_{i}}),t)-\lambda(\psi(t-T_{i},x_{T_{i}}))$.\\
Moreover, conditionally in $(T_{0},x_{T_{0}},\ldots,T_{i},x_{T_{i}},T_{i+1})$, the points in $[T_{i},T_{i+1}[$ obtained from the Poisson process $(T_{i}+\tilde{T}^{i}_{k})_{k\geq{1}}$ are independent of the points in $[T_{j-1},T_{j}[$ obtained from the Poisson process $(T_{j}+\tilde{T}^{j}_{k})_{k\geq{1}}$ for $j=1,\ldots,i$. Thus, we have constructed a point process \[T_{0}<\tilde{T}^{0}_{1}<\ldots< \tilde{T}^{0}_{\tau_{1}-1}< T_{1}<T_{1}+\tilde{T}^{1}_{1}<\ldots<T_{1}+\tilde{T}^{1}_{\tau_{2}-1}<T_{2}<T_{2}+\tilde{T}^{2}_{1}<\ldots\tag{1}\] 
\underline{\textbf{Notation}}: In the sequel, the process defined by $(1)$ is noted $(\tilde{T}_{k})_{k\geq0}$ and the associated counting process is noted $(\tilde{N}_{t})_{t\geq0}$. We also denote by $(\overline{T}_{k})_{k\geq0}$ the process formed by all the rejected points (i.e the process $(\tilde{T}_{k})_{k\geq0}$ without the jump times $(T_{k})_{k\geq0}$) and $(\overline{N}_{t})_{t\geq0}$ the associated counting process.
\vspace{0.2cm}\\
The sequence $(\tilde{T}_{k})_{k\geq0}$ contains both rejected and accepted points and the sub-sequence noted $(T_{k})_{k\geq0}$ such that for $k\geq1$, $T_{k}=\sum_{l=1}^{k}\tilde{T}^{l-1}_{\tau_{l}}$ defines the jump times of the PDMP. Thus, we have constructed the jumps of the PDMP by thinning the process $(\tilde{T}_{k})_{k\geq0}$ with non-constant probabilities $(p_{k})$ such that $p_{k}=\lambda(x_{\tilde{T}_{k}-})/\tilde{\lambda}(x_{\tilde{T}_{k}-},\tilde{T}_{k})$ is the probability to accept $\tilde{T_{k}}$.
\\
An important fact in the procedure is that $(\tilde{T}_{k})_{k\geq0}$ is composed by pieces of independent Poisson processes $(\tilde{T}^{0}_{k}),(\tilde{T}^{1}_{k}),\ldots,(\tilde{T}^{i}_{k}),\ldots$.\\
The thinning procedure provides an algorithm to simulate trajectories of the PDMP.\\
\begin{center}
\textbf{Algorithm}
\end{center}
\rule{\linewidth}{.5pt}
\textbf{Step 1.}\\
Fix the initial time $T_{0}$, the initial condition $x_{T_{0}}=(\theta_{T_{0}},V_{T_{0}})$ and set a jump counter $n=0$. Set also an auxiliary jump counter $k=0$ and an auxiliary variable $\tilde{T}_{k}=T_{n}$.\\
\textbf{Step 2.}\\
$k\leftarrow{k+1}$\\
Simulate $U_{2k-1}\sim\mathcal{U}(]0,1[)$\\
Simulate $E_{k}=-\log(U_{2k-1})$ \\
$\tilde{T}_{k}=\left(\tilde{\Lambda}_{n}\right)^{-1}\left(E_{k}+\tilde{\Lambda}_{n}(\tilde{T}_{k-1})\right)$\\
\textbf{Step 3.}\\
Simulate $U_{2k}\sim\mathcal{U}(]0,1[)$\\
if $U_{2k}\tilde{\lambda}(\psi(\tilde{T}_{k}-T_{n},x_{T_{n}}),\tilde{T}_{k})>\lambda(\psi(\tilde{T}_{k}-T_{n},x_{T_{n}}))$ go to \textbf{Step 2}\\
\textbf{Step 4.}\\
$T_{n+1}=\tilde{T}_{k}$\\
Let $V_{t}=\phi(t-T_{n},x_{T_{n}})$ for $t\in{[T_{n},T_{n+1}[}$. \\
If $T_{n+1}>T$ stop at $t=T$.\\
\textbf{Step 5.}\\
Otherwise, simulate a post-jump value $x_{T_{n+1}}$ according to the Markovian kernel $Q\left(\psi(S_{n+1},x_{T_{n}}),.\right)$.\\
\textbf{Step 6.}\\
Set $n=n+1$ and return to \textbf{Step 2}.\\
\rule{\linewidth}{.5pt}

\section{Efficiency of the thinning procedure}

In this section, we compare the efficiency of the thinning method in term of reject for the different bounds. 
Recall that the number of points needed to simulate one inter-jump time of the PDMP is given by  
\[\tau_{i}^{x}=\inf\{k>0 : U_{k}\tilde{\lambda}^{x}\left(\psi(\tilde{T}^{i-1}_{k},x_{T_{i-1}}),T_{i-1}+\tilde{T}^{i-1}_{k}\right)\leq{\lambda}
\left(\psi(\tilde{T}^{i-1}_{k},x_{T_{i-1}})\right)\}\]
for $i\geq1$, $x\in{\{\text{glo},\text{loc},\text{opt}\}}$ where $(U_{n})_{n\geq1}$ and $(\tilde{T}^{i-1}_{k})_{k\geq0}$ are defined as in section 3.2.
We begin by a lemma.
\vspace{0.3cm}\\
\underline{\textbf{Lemma 4.1}}: We have the following uniform convergence
\[
\sup_{x\in{E}}\sup_{t\geq0}|\tilde{\lambda}^{\epsilon}(t,x)-\lambda\Big{(}\psi(t,x)\Big{)}|\underset{\epsilon\rightarrow0}{\longrightarrow}{0}
\]
where
\[
\tilde{\lambda}^{\epsilon}(t,x)=\sum_{k\geq0}\sup_{s\in{[k\epsilon,(k+1)\epsilon[}}
\lambda\Big{(}\psi(s,x)\Big{)}\textbf{1}_{[k\epsilon,(k+1)\epsilon[}(t)
\]
\begin{proof}
For $n>0$ we set $\epsilon=1/n$, thus \[\tilde{\lambda}^{1/n}(t,x)=\sum_{k\geq0}\sup_{s\in{[k/n,(k+1)/n[}}\lambda\Big{(}\psi(s,x)\Big{)}\textbf{1}_{[k/n,(k+1)/n[}(t).\] 
Let $M=\sup_{x\in{E}}\sup_{t\geq0}\Big{|}\frac{\partial{\lambda}}{\partial_{t}}\Big{(}\psi(t,x)\Big{)}\Big{|}$, $\nu>0$, $N=\lceil{M}/\nu\rceil$ and $n\geq{N}$. Let $x\in{E}$ and $t\geq0$, there exists $l\geq0$ such that $t\in{[l/n,(l+1)/n[}$. Thus 
\[
\tilde{\lambda}^{1/n}(t,x)=\sup_{s\in{[l/n,(l+1)/n[}}\lambda\Big{(}\psi(s,x)\Big{)}.
\]
Let $t_{0}\in{[l/n,(l+1)/n]}$ such that $\sup_{s\in{[l/n,(l+1)/n[}}\lambda\Big{(}\psi(s,x)\Big{)}=
\lambda\Big{(}\psi(t_{0},x)\Big{)}$. The application of the mean value inequality to the function $\lambda\circ\psi$ gives 
\[
|\lambda\Big{(}\psi(t_{0},x)\Big{)}-\lambda\Big{(}\psi(t,x)\Big{)}
|\leq{M}|t_{0}-t|\leq{M}\frac{1}{n}\leq{\nu}.
\]
And the conclusion follows.
\end{proof}

\subsection{Comparison of the mean number of total jump times}

In this section, the $\tau_{i}^{x}$ for $x\in{\{\text{glo},\text{loc},\text{opt}\}}$ are called local reject.
In proposition 4.1, we show that the best local reject is obtained with the \textit{optimal bound}. More the local reject is small less pseudo-random variables have to be simulated. Thus, the simulation time using the \textit{optimal bound} is expected to be smaller than with the two other bounds.
\vspace{0.2cm}\\
\noindent
\underline{\textbf{Proposition 4.1}}: For all $i\geq1$, We have 
\[\mathbb{E}[\tau_{i}^{\text{opt}}|x_{T_{i-1}},T_{i-1}]\leq{\mathbb{E}[\tau_{i}^{\text{loc}}|x_{T_{i-1}},T_{i-1}]}\leq{\mathbb{E}[\tau_{i}^{\text{glo}}|x_{T_{i-1}},T_{i-1}]}.\]
\begin{proof}
Let $i\geq1$ and $S_{i}=T_{i}-T_{i-1}$. We apply theorem 2.2 of \cite{dev} chap.6 and obtain

\begin{align*}
\mathbb{E}[\tau_{i}^{\text{glo}}|x_{T_{i-1}},T_{i-1}]
&=\tilde{\lambda}^{\text{glo}}
\int_{T_{i-1}}^{+\infty}e^{-\int_{T_{i-1}}^{t}\lambda(\psi(s-T_{i-1},x_{T_{i-1}}))ds}dt=
\tilde{\lambda}^{\text{glo}}\mathbb{E}[S_{i}|x_{T_{i-1}},T_{i-1}]\\
\mathbb{E}[\tau_{i}^{\text{loc}}|x_{T_{i-1}},T_{i-1}]
&=\sup_{t\geq{T_{i-1}}}
\lambda\Big{(}\psi(t-T_{i-1},x_{T_{i-1}})\Big{)}
\int_{T_{i-1}}^{+\infty}e^{-\int_{T_{i-1}}^{t}\lambda(\psi(s-T_{i-1},x_{T_{i-1}}))ds}dt\\
&=\sup_{t\geq{T_{i-1}}}
\lambda\Big{(}\psi(t-T_{i-1},x_{T_{i-1}})\Big{)}\mathbb{E}[S_{i}|x_{T_{i}},T_{i-1}]\\
\mathbb{E}[\tau_{i}^{\text{opt}}|x_{T_{i-1}},T_{i-1}]
&=\sum_{k\geq0}\sup_{s\in{\mathcal{P}_{k}^{T_{i-1}}}}
\lambda\Big{(}\psi(s-T_{i-1},x_{T_{i-1}})\Big{)}\int_{\mathcal{P}_{k}^{T_{i-1}}}
e^{-\int_{T_{i-1}}^{t}\lambda(\psi(s-T_{i-1},x_{T_{i-1}}))ds}dt.
\end{align*}
We have
\[
\tilde{\lambda}^{\text{glo}}=\sup_{x\in{E}}
\lambda(x)\geq{\sup_{t\geq{T_{i-1}}}
\lambda\Big{(}\psi(t-T_{i-1},x_{T_{i-1}})\Big{)}}\geq{\sup_{t\in{\mathcal{P}_{k}^{T_{i-1}}}}
\lambda\Big{(}\psi(t-T_{i-1},x_{T_{i-1}})\Big{)}}
\] 
for all $k\geq0$. Since, for all $i\geq1$, $\mathbb{E}[S_{i}|x_{T_{i-1}},T_{i-1}]\geq0$, the conclusion follows.
\end{proof}
\noindent
From Proposition 4.1, we deduce that $\mathbb{E}[\tilde{N}_{t}^{\text{opt}}]\leq{\mathbb{E}[\tilde{N}_{t}^{\text{loc}}]}\leq{\mathbb{E}[\tilde{N}_{t}^{\text{glo}}]}$ where $\tilde{N}_{t}^{\text{opt}}$, $\tilde{N}_{t}^{\text{loc}}$ and $\tilde{N}_{t}^{\text{glo}}$ are counting processes with intensity $\tilde{\lambda}_{t}^{\text{opt}}$, $\tilde{\lambda}_{t}^{\text{loc}}$ and $\tilde{\lambda}_{t}^{\text{glo}}$ respectively.

\subsection{Rate of acceptance}

We are now interested in the rate of acceptance, that is the mean proportion of accepted points in an interval of the form $[0,t]$ for $t>0$. Let $(N_{t})$ be the counting process of the PDMP and $(\tilde{N}_{t})$ the counting process with general jump rate $\tilde{\lambda}$ (see section 3.1). In proposition 4.2 we give an explicit formula of the rate of acceptance defined as $\mathbb{E}[N_{t}/\tilde{N}_{t}|\tilde{N}_{t}\geq1]$. Recall that, for $k\geq1$, 
\[
p_{k}=\lambda(\psi(\tilde{T}_{k}-T_{n_{k}},x_{T_{n_{k}}}))/\tilde{\lambda}(\psi(\tilde{T}_{k}-T_{n_{k}},x_{T_{n_{k}}}),\tilde{T}_{k})
\] 
is the probability to accept the point $\tilde{T}_{k}$ where $T_{n_{k}}$ denotes the last accepted point before $\tilde{T}_{k}$ and $(\tilde{T}_{n})$ is defined by $(1)$. Let $J:\mathbb{R}_{+}\rightarrow\mathbb{R}_{+}$ be the process defined by $J_{t}=\sum_{k\geq0}(t-T_{k})\textbf{1}_{T_{k}\leq{}t<T_{k+1}}$, thus, for $t\geq0$, $J_{t}$ gives the age of the last accepted point before $t$. Then, for $k\geq1$, we can write the probabilities $p_{k}$ as follows    
\[
p_{k}=\lambda(\psi(J_{\tilde{T}_{k-1}}+\tilde{S_{k}},\eta_{\tilde{T}_{k-1}}))/\tilde{\lambda}(\psi(J_{\tilde{T}_{k-1}}+\tilde{S_{k}},\eta_{\tilde{T}_{k-1}}),\tilde{T}_{k})
\]
where $\tilde{S_{k}}=\tilde{T}_{k}-\tilde{T}_{k-1}$ and $(\eta_{t})$ is as in section 2. The process $(\tilde{S_{k}},\tilde{T}_{k},J_{\tilde{T}_{k}},\eta_{\tilde{T}_{k}})_{k\geq0}$ defines a Markov chain on $\mathbb{R}_{+}\times\tilde{E}$ where $\tilde{E}=\mathbb{R}_{+}\times\mathbb{R}_{+}\times{E}$ with semi-Markovian kernel 
\[M(\tilde{T}_{k-1},J_{\tilde{T}_{k-1}},\eta_{\tilde{T}_{k-1}},ds,dt,dj,dx)=\tilde{Q}(s,\tilde{T}_{k-1},J_{\tilde{T}_{k-1}},\eta_{\tilde{T}_{k-1}},dt,dj,dx)\alpha(\tilde{T}_{k-1},J_{\tilde{T}_{k-1}},\eta_{\tilde{T}_{k-1}},ds)\]
where 
\[
\alpha(\tilde{T}_{k-1},J_{\tilde{T}_{k-1}},\eta_{\tilde{T}_{k-1}},ds)=\tilde{\lambda}(\psi(J_{\tilde{T}_{k-1}}+s,\eta_{\tilde{T}_{k-1}}),\tilde{T}_{k-1}+s)e^{-\int_{0}^{s}\tilde{\lambda}(\psi(J_{\tilde{T}_{k-1}}+z,\eta_{\tilde{T}_{k-1}}),\tilde{T}_{k-1}+z)dz}ds
\]
and 
\begin{align*}
&\tilde{Q}(s,\tilde{T}_{k-1},J_{\tilde{T}_{k-1}},\eta_{\tilde{T}_{k-1}},dt,dj,dx)=\\
&\Big{(}1-\frac{\lambda(\psi(J_{\tilde{T}_{k-1}}+s,\eta_{\tilde{T}_{k-1}}))}{\tilde{\lambda}(\psi(J_{\tilde{T}_{k-1}}+s,\eta_{\tilde{T}_{k-1}}),t)}\Big{)}\delta_{J_{\tilde{T}_{k-1}}+s}(dj)\delta_{\eta_{\tilde{T}_{k-1}}}(dx)\delta_{\tilde{T}_{k-1}+s}(dt)+\\
&\frac{\lambda(\psi(J_{\tilde{T}_{k-1}}+s,\eta_{\tilde{T}_{k-1}}))}{\tilde{\lambda}(\psi(J_{\tilde{T}_{k-1}}+s,\eta_{\tilde{T}_{k-1}}),t)}Q\Big{(}\lambda(\psi(J_{\tilde{T}_{k-1}}+s,\eta_{\tilde{T}_{k-1}})),dx\Big{)}\delta_{0}(dj)\delta_{\tilde{T}_{k-1}+s}(dt)
\end{align*}
\vspace{0.2cm}\\
\textbf{\underline{Proposition 4.2}}: Let $(N_{t})_{t\geq0}$ be the counting process of the PDMP $(x_{t})_{t\geq0}$, $(\tilde{N}_{t})_{t\geq0}$ be the counting process with jump times $(\tilde{T}_{n})_{n\geq0}$ and $M$ be the kernel of the Markov chain $(\tilde{S_{k}},\tilde{T_{k}},J_{\tilde{T}_{k}},\eta_{\tilde{T}_{k}})_{k\geq0}$, we have
\begin{align*}
\mathbb{E}\Big{[}\frac{N_{t}}{\tilde{N}_{t}}|\tilde{N}_{t}\geq1\Big{]}=\frac{1}{\mathbb{P}(\tilde{N}_{t}\geq1)}\sum_{n\geq1}\frac{1}{n}\int_{\mathbb{R}_{+}}\int_{\tilde{E}}\sum_{k=1}^{n}\left(
\frac{\lambda\Big{(}\psi(j_{k-1}+s_{k},x_{k-1})\Big{)}}{\tilde{\lambda}\Big{(}\psi(j_{k-1}+s_{k},x_{k-1}),t_{k}\Big{)}}\right)e^{-\int_{0}^{t-t_{n}}\tilde{\lambda}(\psi(j_{n}+z,x_{n}),t_{n}+z)dz}
&\\\textbf{1}_{t\geq{t_{n}}}\mu(dx_{0})
M(0,0,x_{0},ds_{1},dt_{1},dj_{1},dx_{1})\ldots{M}(t_{n-1},j_{n-1},x_{n-1},ds_{n},dt_{n},dj_{n},dx_{n})
\end{align*}
where $\mu$ is the law of $\eta_{\tilde{T}_{0}}$.

\begin{proof}

We provide a proof into two steps. First, we establish that, with an appropriate conditioning, the conditional law of $N_{t}$ is the conditional law of a sum of independent Bernoulli random variables with different parameters. Then, we use this property as well as the kernel $M$ to compute the rate of acceptance.\\
Let $n\geq1$ and let us define $n$ independent Bernoulli random variables $X_{i}$ with parameters $p_{i}$ such that
\[
p_{i}=\frac{\lambda(\psi(J_{\tilde{T}_{i-1}}+\tilde{S}_{i},\eta_{\tilde{T}_{i-1}}))}{\tilde{\lambda}(\psi(J_{\tilde{T}_{i-1}}+\tilde{S}_{i},\eta_{\tilde{T}_{i-1}}),\tilde{T}_{i})}
\]
We note $X=\sum_{i=1}^{n}X_{i}$ the sum of the Bernoulli variables and $A_{t,n}=\{\tilde{N}_{t}=n,p_{1},\ldots,p_{n}\}$. 
By noting that, for $0\leq{k}\leq{n}$, we have
\[
\{N_{t}=k|A_{t,n}\}
=\bigcup_{1\leq{i_{1}}<\ldots<i_{k}\leq{n}}\Big{[}
\bigcap_{i\in{I_{k}^{n}}}\{U_{i}\leq{p_{i}}\}
\bigcap_{i\in{\overline{I}_{k}^{n}}}\{U_{i}>p_{i}\}\Big{]}
=\{X=k|p_{1},\ldots,p_{n}\}
\]
where $I_{k}^{n}=\{i_{1},\ldots,i_{k}\}\subseteq{\{1,\ldots,n\}}$, $\overline{I}_{k}^{n}$ the complementary of $I_{k}^{n}$ in $\{1,\ldots,n\}$ and $(U_{i})$ are independent random variables uniformly distributed in $[0,1]$ and independent of $(p_{i})$, we deduce that
\[
\mathcal{L}(N_{t}|A_{t,n})=\mathcal{L}(X|p_{1},\ldots,p_{n}).
\]
In particular, $\mathbb{E}[N_{t}|A_{t,n}]=\mathbb{E}[X|p_{1},\ldots,p_{n}]=\sum_{i=1}^{n}p_{i}$.
Thus, one can write 
\begin{align*}
\mathbb{E}\Big{[}\frac{N_{t}}{\tilde{N}_{t}}|\tilde{N}_{t}\geq1\Big{]}
&=\frac{1}{\mathbb{P}(\tilde{N}_{t}\geq1)}\sum_{n\geq1}\frac{1}{n}\mathbb{E}\Big{[}N_{t}\textbf{1}_{\tilde{N}_{t}=n}\Big{]}\\
&=\frac{1}{\mathbb{P}(\tilde{N}_{t}\geq1)}\sum_{n\geq1}\frac{1}{n}\mathbb{E}\Big{[}\mathbb{E}\Big{[}N_{t}|A_{t,n}\Big{]}\Big{|}\tilde{N}_{t}=n\Big{]}\mathbb{P}(\tilde{N}_{t}=n)\\
&=\frac{1}{\mathbb{P}(\tilde{N}_{t}\geq1)}\sum_{n\geq1}\frac{1}{n}\mathbb{E}\Big{[}\sum_{i=1}^{n}p_{i}\textbf{1}_
{\tilde{N}_{t}=n}\Big{]}\\
&=\frac{1}{\mathbb{P}(\tilde{N}_{t}\geq1)}\sum_{n\geq1}\frac{1}{n}\mathbb{E}\Big{[}\sum_{i=1}^{n}p_{i}\textbf{1}_
{t-\tilde{T}_{n}\geq{0}}\mathbb{E}\Big{[}\textbf{1}_
{\tilde{S}_{n+1}\geq{t-\tilde{T}_{n}}}|\eta_{\tilde{T}_{0}},\tilde{S_{1}},\ldots,\tilde{S_{n}},\tilde{T_{n}},J_{\tilde{T}_{n}},\eta_{\tilde{T}_{n}}\Big{]}\Big{]}\\
\intertext{Conditionally to $(\eta_{\tilde{T}_{0}},\tilde{S_{1}},\tilde{T_{1}},J_{\tilde{T}_{1}},\eta_{\tilde{T}_{1}},\ldots,\tilde{S_{n}},\tilde{T_{n}},J_{\tilde{T}_{n}},\eta_{\tilde{T}_{n}})$, the random variable $\tilde{S}_{n+1}$ is a hazard law with rate $\tilde{\lambda}(\psi(J_{\tilde{T}_{n}}+t,\eta_{\tilde{T}_{n}}),\tilde{T}_{n}+t)$ for $t\geq{0}$, thus}
\mathbb{E}\Big{[}\frac{N_{t}}{\tilde{N}_{t}}|\tilde{N}_{t}\geq1\Big{]}&=\frac{1}{\mathbb{P}(\tilde{N}_{t}\geq1)}\sum_{n\geq1}\frac{1}{n}\mathbb{E}\Big{[}\sum_{i=1}^{n}p_{i}e^{-\int_{0}^{t-\tilde{T}_{n}}\tilde{\lambda}(\psi(J_{\tilde{T}_{n}}+u,\eta_{\tilde{T}_{n}}),\tilde{T}_{n}+u)du}
\textbf{1}_{t-\tilde{T}_{n}\geq0}\Big{]}\\
&=\frac{1}{\mathbb{P}(\tilde{N}_{t}\geq1)}\sum_{n\geq1}\frac{1}{n}\mathbb{E}\Big{[}f(\eta_{\tilde{T}_{0}},\tilde{S_{1}},\tilde{T_{1}},J_{\tilde{T}_{1}},\eta_{\tilde{T}_{1}},\ldots,\tilde{S_{n}},\tilde{T}_{n},J_{\tilde{T}_{n}},\eta_{\tilde{T}_{n}})\Big{]}
\end{align*}
Where
\vspace{0.1cm}\\
$f(\eta_{\tilde{T}_{0}},\tilde{S_{1}},\tilde{T_{1}},J_{\tilde{T}_{1}},\eta_{\tilde{T}_{1}},\ldots,\tilde{S_{n}},\tilde{T_{n}},J_{\tilde{T}_{n}},\eta_{\tilde{T}_{n}})=$
\[=\sum_{i=1}^{n}\left(\frac{\lambda(\psi(J_{\tilde{T}_{i-1}}+\tilde{S}_{i},\eta_{\tilde{T}_{i-1}}))}{\tilde{\lambda}(\psi(J_{\tilde{T}_{i-1}}+\tilde{S}_{i},\eta_{\tilde{T}_{i-1}}),\tilde{T}_{i})}\right)e^{-\int_{0}^{t-\tilde{T}_{n}}\tilde{\lambda}(\psi(J_{\tilde{T}_{n}}+u,\eta_{\tilde{T}_{n}}),\tilde{T}_{n}+u)du}\textbf{1}_{t-\tilde{T}_{n}\geq0}\]
Since $(\tilde{S_{k}},\tilde{T_{k}},J_{\tilde{T}_{k}},\eta_{\tilde{T}_{k}})_{k\geq0}$ is a Markov chain with kernel $M$, we obtain
\begin{align*}
\mathbb{E}\left[f(\eta_{\tilde{T}_{0}},\tilde{S_{1}},\ldots,\tilde{S_{n}},J_{\tilde{T}_{n}},\eta_{\tilde{T}_{n}})\right]=\int_{\mathbb{R}_{+}}\int_{\tilde{E}}
f(x_{0},s_{1},t_{1},j_{1},x_{1}\ldots,s_{n},t_{n},j_{n},x_{n})
&\\\mu(dx_{0})
M(0,0,x_{0},ds_{1},dt_{1},dj_{1},dx_{1})\ldots
M(t_{n-1},j_{n-1},x_{n-1},ds_{n},dt_{n},dj_{n},dx_{n})
\end{align*}
Where $\mu$ is the law of $\eta_{\tilde{T}_{0}}$. Thus we have the result.
\end{proof}
\noindent
When $\tilde{\lambda}$ is close to $\lambda$, the rate of acceptance is expected to be close to 1. As an example, consider the case of two Poisson processes $(N_{t})$ and $(\tilde{N}_{t})$ with intensity $\lambda(t)$ and $\tilde{\lambda}(t)$ respectively such that $\tilde{\lambda}(t)=\tilde{\lambda}$ for all $t\geq0$. Thus, for $n\geq1$, $\tilde{S}_{1},\ldots,\tilde{S}_{n}$ are independent exponential variables with parameter $\tilde{\lambda}$. Let us also consider that $\lambda(t)\simeq\tilde{\lambda}$ for $t\geq0$. In this case, the rate of acceptance is
\[
\mathbb{E}\Big{[}\frac{N_{t}}{\tilde{N}_{t}}|\tilde{N}_{t}\geq1\Big{]}\simeq\frac{1}{1-e^{-\tilde{\lambda}t}}\sum_{n\geq1}\int_{(\mathbb{R}_{+})^{n}}
e^{-\tilde{\lambda}\Big{(}t-(s_{1}+\ldots+s_{n})\Big{)}}
\textbf{1}_{t\geq{s_{1}+\ldots+s_{n}}}\alpha(ds_{1}),\ldots,\alpha(ds_{n})
\]
where $\alpha(ds)=\tilde{\lambda}e^{-\tilde{\lambda}s}ds$. Since, $\tilde{T}_{n}=\tilde{S}_{1}+\ldots+\tilde{S}_{n}$ is gamma distributed with parameters $n$ and $\tilde{\lambda}$, we have
\[
\mathbb{E}\Big{[}\frac{N_{t}}{\tilde{N}_{t}}|\tilde{N}_{t}\geq1\Big{]}\simeq\frac{1}{1-e^{-\tilde{\lambda}t}}\sum_{n\geq1}\mathbb{E}[e^{-\tilde{\lambda}(t-{\tilde{T}_{n}})}\textbf{1}_{t\geq\tilde{T}_{n}}]
\simeq\frac{1}{1-e^{-\tilde{\lambda}t}}\sum_{n\geq1}\frac{(\tilde{\lambda}t)^{n}}{n!}e^{-\tilde{\lambda}t}
\simeq1
\]

\subsection{Convergence of the counting process with a specific optimal bound as jump rate}

Let $P=\mathbb{N}$ and let $\epsilon>0$, we consider the optimal bound with the partition $(\mathcal{P}_{k}^{T_{n},\epsilon})_{k\in{\mathbb{N}}}$ where, for $k\in{\mathbb{N}}$, $\mathcal{P}_{k}^{T_{n},\epsilon}=[T_{n}+k\epsilon,T_{n}+(k+1)\epsilon[$, we denote it by $\tilde{\lambda}^{\text{opt},\epsilon}(x_{t},t)$ and we note $(\tilde{N}^{\text{opt},\epsilon}_{t})$ the corresponding counting process. The number of points needed to simulate one inter-jump time with this particular partition is noted $\tau_{i}^{\text{opt},\epsilon}$.\\
We first show, in proposition 4.3, that $(\overline{T}_{k})_{k\geq0}$ (defined in section 3.2) is a Cox process with stochastic jump rate $\tilde{\lambda}(x_{t},t)-\lambda(x_{t})$. Details on Cox processes can be found in \cite{kal}.

Then we show, in proposition 4.4, that the counting process $(\tilde{N}^{\epsilon}_{t})$ whose jump rate is $\tilde{\lambda}^{\text{opt},\epsilon}$ converge to the counting process $(N_{t})$ of the PDMP. 

Proposition 4.5 states that more the parameter $\epsilon$ is small less rejected points are simulated. However, when $\epsilon$ is too small, step 2 of the algorithm requires many iterations to compute $\tilde{\Lambda}(.)$ and $\left(\tilde{\Lambda}\right)^{-1}(.)$ and the simulation time increases. We will see in the numerical section 6.2.1 that taking an $\epsilon$ of order $\max_{n}(T_{n+1}-T_{n})$ leads to the optimal simulation time. 
\vspace{0.2cm}\\
\underline{\textbf{Proposition 4.3}} : Let $\xi$ be a point process and $\mu$ be a random measure such that 
\[
\xi([0,t])=\sum_{n\geq0}\textbf{1}_{\overline{T}_{n}
\leq{t}}\hspace{0.3cm}\text{and}\hspace{0.3cm}
\mu([0,t])=\int_{0}^{t}\tilde{\lambda}(x_{s},s)-\lambda(x_{s})ds
\]
Then $\xi$ is a Cox process directed by $\mu$.
\begin{proof}
Let us first note that for $t\geq0$
\[\tilde{\lambda}(x_{t},t)-\lambda(x_{t})=\sum_{n\geq0}\Big{[}\tilde{\lambda}(\psi(t-T_{n},x_{T_{n}}),t)-\lambda(\psi(t-T_{n},x_{T_{n}}))\Big{]}\textbf{1}_{T_{n}\leq{t}<T_{n+1}}\]
Thanks to \cite{kal}, we show that for any measurable and non-negative functions $f$, the Laplace functional of $\xi$ is
\[\mathbb{E}[e^{-\xi{f}}]=\mathbb{E}[e^{-\mu(1-e^{-f})}]
\] 
Let $f$ be a non-negative measurable function. Let us note $f_{T}(t)=f(t)\textbf{1}_{t\leq{T}}$ for $T>0$, such that $\lim_{T\rightarrow\infty}f_{T}(t)=f(t)$ with $f_{T}$ increasing with $T$. Thus, by Beppo-Levi's theorem, $\xi{f}_{T}\nearrow\xi{f}$ and then $e^{-\xi{f}_{T}}\searrow{e}^{-\xi{f}}$ when $T$ goes to infinity.
Moreover, $e^{-\xi{f}_{T}}\leq1$, thus by Lebesgue's Dominated Convergence Theorem
\[
\mathbb{E}[e^{-\xi{f}_{T}}]\rightarrow\mathbb{E}[e^{-\xi{f}}]
\]
With the same type of arguments, we show that
\[\mathbb{E}[e^{-\mu(1-e^{-f_{T}})}]\rightarrow\mathbb{E}[e^{-\mu(1-e^{-f})}]
\]
Thus, it is sufficient to show $(2)$ for functions $f_{T}$.
\begin{align*}
\mathbb{E}[e^{-\xi{f}_{T}}]&=\mathbb{E}[e^{-\sum_{n\geq1}f_{T}(\overline{T}_{n})}]\\
&=\sum_{k\geq0}\mathbb{E}[e^{-\sum_{n\geq1}f_{T}(\overline{T}_{n})}|N_{T}=k]\mathbb{P}(N_{T}=k)\\
&=\sum_{k\geq0}\mathbb{E}\Big{[}\mathbb{E}[\prod_{i=0}^{k}e^{-\sum_{n\geq1}f_{T}(\overline{T}_{n})\textbf{1}_{T_{i}\leq{\overline{T}_{n}}<T_{i+1}}}|N_{T}=k,(\eta_{t})_{0\leq{t}\leq{T}}]|N_{T}=k\Big{]}\mathbb{P}(N_{T}=k)
\end{align*}
By the thinning procedure, the points $\overline{T}_{n}$ in $[T_{i},T_{i+1}[$ may be written as $T_{i}+\tilde{T}^{i}_{l}$ for some $l\geq1$ where $(\tilde{T}^{i}_{l})_{l\geq1}$ is, conditionally on $x_{T_{i}}$, a Poisson process with jump rate $\tilde{\lambda}(\psi(t-T_{i},x_{T_{i}}),t)-\lambda(\psi(t-T_{i},x_{T_{i}}))$ for $t\geq{T_{i}}$. Since $(\tilde{T}^{i}_{l})$ is independent of $(\tilde{T}^{j}_{l})$ for $i\neq{j}$, the random variables $X_{i}:=e^{-\sum_{n\geq1}f_{T}(\overline{T}_{n})\textbf{1}_{T_{i}\leq{\overline{T}_{n}}<T_{i+1}}}$ are independent. Thus, since the Laplace functional of a Poisson process $\xi$ with intensity $\mu$ verifies $\mathbb{E}[e^{-\xi{f}}]=e^{-\mu(1-e^{-f})}$, we obtain

\begin{align*}
\mathbb{E}[e^{-\xi{f}_{T}}]&=\sum_{k\geq0}\mathbb{E}\Big{[}\prod_{i=0}^{k}\mathbb{E}[e^{-\sum_{n\geq1}f_{T}(\overline{T}_{n})\textbf{1}_{T_{i}\leq{\overline{T}_{n}}<T_{i+1}}}|N_{T}=k,(\eta_{t})_{0\leq{t}\leq{T}}]|N_{T}=k\Big{]}\mathbb{P}(N_{T}=k)\\
&=\sum_{k\geq0}\mathbb{E}\Big{[}
e^{-\sum_{i=0}^{N_{T}}\int\Big{(}1-e^{-f_{T}(s)\textbf{1}_{T_{i}\leq{s}<T_{i+1}}}\Big{)}\Big{(}\tilde{\lambda}(\psi(s-T_{i},x_{T_{i}}),s)-\lambda(\psi(s-T_{i},x_{T_{i}}))\Big{)}ds}|N_{T}=k\Big{]}\mathbb{P}(N_{T}=k)\\					
&=\mathbb{E}\Big{[}e^{-\sum_{i\geq0}\int\Big{(}1-e^{-f_{T}(s)}\Big{)}\textbf{1}_{T_{i}\leq{}s<T_{i+1}}\Big{(}\tilde{\lambda}(\psi(s-T_{i},x_{T_{i}}),s)-\lambda(\psi(s-T_{i},x_{T_{i}}))\Big{)}ds}\Big{]}\\						
&=\mathbb{E}[e^{-\mu(1-e^{-f_{T}})}]\\
\end{align*}
\end{proof}

\noindent
\underline{\textbf{Proposition 4.4}}: Let $(\tilde{N}^{\epsilon}_{t})$ be the counting process whose jump rate is $\tilde{\lambda}^{\text{opt},\epsilon}$ and let $(N_{t})$ be the counting process of the PDMP $x_{t}$, we have the following convergence in law
\[
\tilde{N}^{\epsilon}\underset{\epsilon\rightarrow0}{\longrightarrow}{N}
\]
\begin{proof}
Let $f$ be a non-negative measurable function . We show the convergence of the Laplace transform \cite{cocofiab}, that is 
\[
\mathbb{E}[e^{-\int{f}d\tilde{N}^{\epsilon}}]\underset{\epsilon\rightarrow0}{\longrightarrow}\mathbb{E}[e^{-\int{f}dN}].
\]
Let $T>0$, as in proposition 4.3, it is sufficient to show the convergence of the Laplace transform for functions $f_{T}(t)=f(t)\textbf{1}_{t\leq{T}}$. Let $(\tilde{T}^{\epsilon}_{n})$ be the points of the process $\tilde{N}^{\epsilon}$.
\begin{align*}
\mathbb{E}[e^{-\int{f_{T}}d\tilde{N}^{\epsilon}}]&=\mathbb{E}[e^{-\sum_{n\geq0}f_{T}(\tilde{T}^{\epsilon}_{n})}]\\
&=\mathbb{E}\Big{[}\mathbb{E}[e^{-\sum_{n\geq0}f_{T}(\tilde{T}^{\epsilon}_{n})}|(\eta_{t})_{0\leq{t}\leq{T}}]\Big{]}\\
&=\mathbb{E}\Big{[}e^{-\sum_{n\geq0}f_{T}(T_{n})}\mathbb{E}[e^{-
\sum_{n\geq0}f_{T}(\overline{T}^{\epsilon}_{n})
}|(\eta_{t})_{0\leq{t}\leq{T}}]\Big{]}\\
\intertext{Where $(\overline{T}_{n}^{\epsilon})$ denotes the rejected points. Since $(\overline{T}_{n}^{\epsilon})$ is a Cox process with stochastic jump rate $\tilde{\lambda}^{\text{opt},\epsilon}(x_{t},t)-\lambda(x_{t})$, we obtain}
\mathbb{E}[e^{-\int{f_{T}}d\tilde{N}^{\epsilon}}]&=\mathbb{E}\Big{[}e^{-\sum_{n\geq0}f_{T}(T_{n})}
e^{-\int(1-e^{-f_{T}(s)})(\tilde{\lambda}^{\text{opt},\epsilon}(x_{s},s)-\lambda(x_{s}))ds}\Big{]}
\end{align*}
Since $e^{-\sum_{n\geq0}f_{T}(T_{n})}
e^{-\int(1-e^{-f_{T}(s)})(\tilde{\lambda}^{\text{opt},\epsilon}(x_{s},s)-\lambda(x_{s}))ds}\leq1$, we obtain by Lebesgue's Dominated Convergence Theorem and by continuity of the exponential
\[
\lim_{\epsilon\rightarrow0}\mathbb{E}[e^{-\int{f_{T}}d\tilde{N}^{\epsilon}}]=\mathbb{E}\Big{[}e^{-\sum_{n\geq0}f_{T}(T_{n})}
e^{\lim_{\epsilon\rightarrow0}\int-(1-e^{-f_{T}(s)})(\tilde{\lambda}^{\text{opt},\epsilon}(x_{s},s)-\lambda(x_{s}))ds}\Big{]}.
\]
We have
\[
-T\sup_{x\in{E}}\sup_{t\geq0}\Big{(}\tilde{\lambda}^{\epsilon}(t,x)-\lambda(\psi(t,x))\Big{)}\leq{}\int-(1-e^{-f_{T}(s)})(\tilde{\lambda}^{\text{opt},\epsilon}(x_{s},s)-\lambda(x_{s}))ds\leq0
\]
where 
\[
\tilde{\lambda}^{\epsilon}(t,x)=\sum_{k\geq0}\sup_{s\in{[k\epsilon,(k+1)\epsilon[}}
\lambda\Big{(}\psi(s,x)\Big{)}\textbf{1}_{[k\epsilon,(k+1)\epsilon[}(t)
\]
By lemma 4.1, we obtain that almost surely $e^{\lim_{\epsilon\rightarrow0}\int-(1-e^{-f_{T}(s)})(\tilde{\lambda}^{\text{opt},\epsilon}(x_{s},s)-\lambda(x_{s}))ds}=1$, the conclusion follows since $\mathbb{E}\Big{[}e^{-\sum_{n\geq0}f_{T}(T_{n})}\Big{]}=\mathbb{E}\Big{[}e^{-\int{f}_{T}dN}\Big{]}$.
\end{proof}
\noindent
\underline{\textbf{Proposition 4.5}}: For all $i\geq1$ we have 
\[
\mathbb{E}[\tau_{i}^{\text{opt},\epsilon}|x_{T_{i-1}},T_{i-1}]\underset{\epsilon\rightarrow0}{\longrightarrow}{1}
\]
\begin{proof}
Let $i\geq1$ and $\epsilon>0$.
From theorem 2.2 in chap.6 of \cite{dev} we have 
\[
\mathbb{E}[\tau_{i}^{\text{opt},\epsilon}|x_{T_{i-1}},T_{i-1}]=\int_{T_{i-1}}^{+\infty}
\tilde{\lambda}^{\epsilon}_{i-1}(t-T_{i-1})e^{-\int_{T_{i-1}}^{t}\lambda(\psi(s-T_{i-1},x_{T_{i-1}}))ds}dt
\]
Where
\[
\tilde{\lambda}^{\epsilon}_{i-1}(t-T_{i-1})=\sum_{k\geq0}\sup_{s\in{\mathcal{P}_{k}^{T_{i-1}}}}
\lambda\Big{(}\psi(s-T_{i-1},x_{T_{i-1}})\Big{)}\textbf{1}_{\mathcal{P}_{k}^{T_{i-1},\epsilon}}(t)
\]
from theorem 2.3 in \cite{dev} chap.6, we deduce that
\[
\mathbb{E}[\tau_{i}^{\text{opt},\epsilon}|x_{T_{i-1}},T_{i-1}]\leq{}\sup_{t\geq{T_{i-1}}}
\frac{\tilde{\lambda}^{\epsilon}_{i-1}(t-T_{i-1})}{\lambda\Big{(}\psi(t-T_{i-1},x_{T_{i-1}})\Big{)}}. 
\]
By Lemma 4.1, we obtain 
\[
\lim_{\epsilon\rightarrow0}\mathbb{E}[\tau_{i}^{\text{opt},\epsilon}|x_{T_{i-1}},T_{i-1}]\leq{1}.
\]
Since $\mathbb{E}[\tau_{i}^{\text{opt},\epsilon}|x_{T_{i-1}},T_{i-1}]\geq1$ for all $\epsilon>0$, the conclusion follows.
\end{proof}

\section{Hodgkin-Huxley models}

In this section we apply the thinning algorithms above to different stochastic versions of the Hodgkin-Huxley model. First, we introduce the classical deterministic model and the stochastic models, then we give essential points of the simulation.

\subsection{Deterministic Hodgkin-Huxley model}

Alan Lloyd Hodgkin and Andrew Huxley provided two models in $1952$ to explain the ionic mechanisms underlying the initiation and propagation of action potentials in the squid giant axon \cite{hodgkin}.\\
In the present section, we only consider models of initiation of action potentials. It means that we clamp (isolate) a piece of the axon or of the soma and we study electrical properties in time in this clamped area (also called membrane-patch).\\
Neuron's membrane separates the intracellular environment from the extracellular one and allows exchanges of material and energy between these two environments. These exchanges are allowed by the opening and closing of gates located on the membrane.  
In most neurons, the intracellular environment contains a large proportion of potassium ions, whereas the extracellular environment contains a majority of sodium ones.
Hodgkin and Huxley discovered that the initiation of action potentials principally relies on the movement of these two kind of ions across the membrane via ionic channels. A ionic channel is constituted by four gates which can be of different types (activation and inactivation) and is specific to one type of ions, for example, a sodium channel allows sodium ions only to pass the membrane. We say that a channel is active when all his gates are open.\\ 
A stimulation (it can be an input from other neurons or external applied current) makes the sodium channels active, thus sodium ions enter in the intracellular environment : the membrane is depolarized. This depolarizing step increases the conductance of the membrane and when the voltage exceeds a threshold value there is an action potential. After being active, sodium channels become inactive, while potassium gates open (these opening make the potassium channels active). Potassium ions leave the intracellular environment to compensate the entry of sodium ions : the membrane is re-polarized. Potassium channels stay active longer than sodium ones : the membrane is hyper-polarized. Then, a protein makes the sodium ions go back into the intracellular environment and expels potassium ions outside. These are the principal steps of the initiation of an action potential.

\begin{figure}
\includegraphics[width=0.55\textwidth]{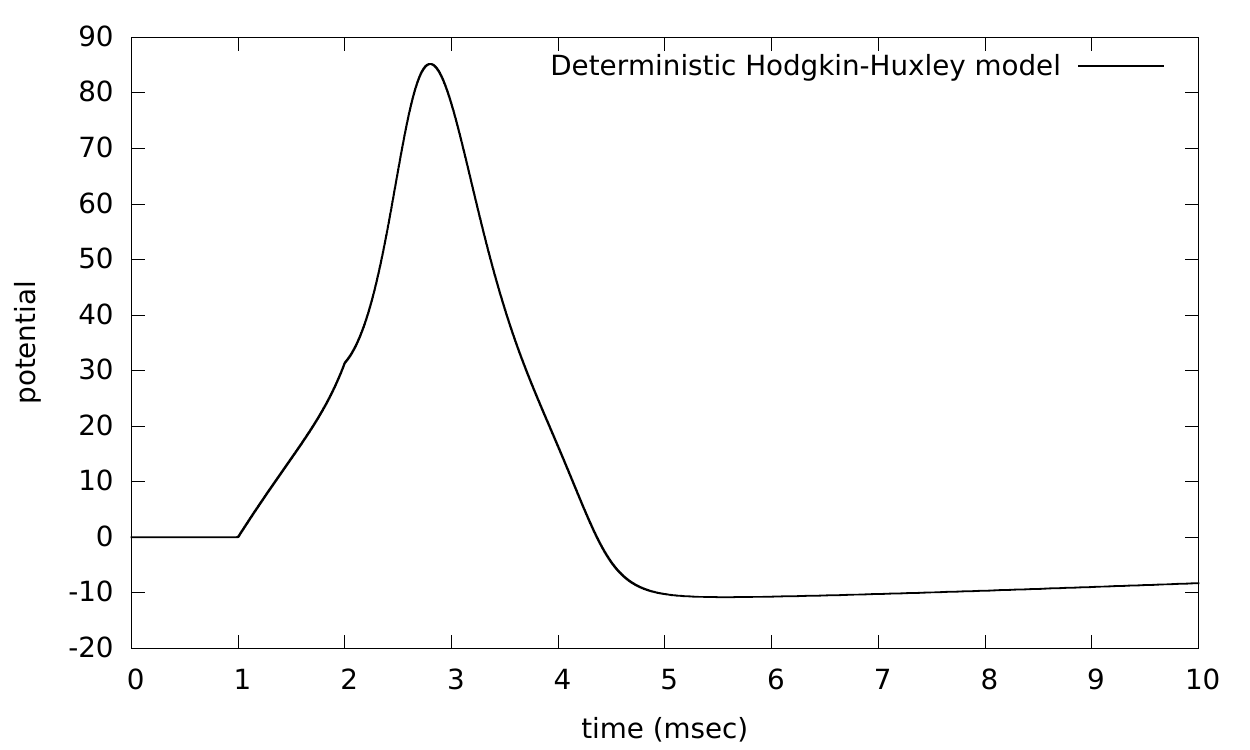}
\caption{A trajectory of the deterministic Hodgkin-Huxley model with the values of the parameters given in appendix A and $I(t)=30\textbf{1}_{[1,2]}(t)$.}
\label{fig:1}
\end{figure}

In the models studied in this section, we distinguish three types of channels : sodium, potassium, and leak. We consider that a sodium (potassium) channel is composed by three activation gates represented by the variable $m$ and one inactivation gate $h$ (four activation gates $n$ and zero inactivation gates) and that leak channels are always open and allow all type of ions to pass the membrane. \\
In fact, there exists other ionic channels as for example calcium ones. In some cases, the calcium plays the role of the sodium to initiate action potentials as in the giant barnacle muscle fiber \cite{lecar}, or more generally in crustacean muscles \cite{Bertil}. Calcium channels have also many other functionalities \cite{Bertil}.
Classically, the Hodgkin-Huxley model is the set of nonlinear differential equations
\begin{align*}
C&\frac{dV}{dt}=I-I_{L}(V)-I_{Na}(V,m,h)-I_{K}(V,n)\\
&\frac{dm}{dt}=(1-m)\alpha_{m}(V)-m\beta_{m}(V)\\
&\frac{dh}{dt}=(1-h)\alpha_{h}(V)-h\beta_{h}(V)\\
&\frac{dn}{dt}=(1-n)\alpha_{n}(V)-n\beta_{n}(V)
\end{align*}
The function $V$ represents the membrane potential (voltage). The functions $m$, $h$, $n$ correspond to the fraction of open gates of type $m$, $h$, or $n$. The functions $\alpha_{x}$ and $\beta_{x}$ for $x=m,h,n$ are opening and closing rates of gates $x$ respectively. $I$ is a time-dependent function which represents the input current, $C$ is the membrane capacity. For $z\in\{\text{Na},\text{K},\text{L}\}$, $I_{z}=\overline{g}_{z}(V-V_{z})$ represents the ionic currents where $\overline{g}_{\text{Na}}=g_{\text{Na}}m^{3}h$, $\overline{g}_{\text{K}}=g_{\text{K}}n^{4}$ and $\overline{g}_{\text{L}}=g_{\text{L}}$ are the conductances of the sodium, potassium and leak respectively. $g_{\text{L}}$, $g_{\text{Na}}$, $g_{\text{K}}$ are the conductances when all gates are opened and $V_{\text{L}}$, $V_{\text{Na}}$, $V_{\text{K}}$ are the resting potentials. Note that, the conductance of the membrane depends on the state of the gates. Thus, more the open gates are, more the conductance is high.
This model describes the electrical behaviour of a neuron with an infinite number of gates. Thus, it does not reflect the variability observed experimentally.

\subsection{Stochastic Hodgkin-Huxley models}

Neurons are subject to various sources of fluctuations, intrinsic (from the membrane) and extrinsic (from synapses). The intrinsic fluctuations are mainly caused by ion channels. To take into account these fluctuations in the model, we fix a finite number of sodium channels $N_{\text{Na}}$ and potassium ones $N_{\text{K}}$, and replace the deterministic dynamic of the gates by stochastic processes. 
Here, we discuss about two stochastic models, namely the $\textit{subunit model}$ and the $\textit{channel model}$. These models belong to the class of Piecewise Deterministic Markov Processes.

\subsubsection{The subunit model}

The \textit{subunit model} is obtained by considering that the conductance of the membrane depends on the empirical measure defined by the proportion of open gates. 
Recall that $m$ (respectively $n$) denotes an activation gate of the sodium (respectively potassium) channels and $h$ an inactivation one. We define the number of gates $m$ (respectively $h$, $n$) by $N_{m}=3N_{\text{Na}}$ (respectively $N_{h}=N_{\text{Na}}$, $N_{n}=4N_{\text{K}}$). Let us consider that each gate is represented by a $\{0,1\}$-valued Markovian Jump Process (MJP) noted $u_{k}^{(x)}$ for $x=m,h,n$ and $k=1,\ldots,N_{x}$. State 1 corresponds to the open configuration and 0 to the close one. Recall also that the opening and closing rates which depend on the voltage are noted $\alpha_{x}(.)$ and $\beta_{x}(.)$ respectively. The dynamic of a gate can be represented by the following diagram
\[
\begin{array}{lll}
0\begin{array}{l} \alpha_{x}(.)\\\longrightarrow\\\longleftarrow \\ \beta_{x}(.) \end{array}1
\end{array}
\]
We consider that all MJPs are independent and we define the number of open gates $x$ at time $t$ by 
\[
\theta^{(x)}(t)=\sum_{k=1}^{N_{x}}u_{k}^{(x)}(t).
\] 
Furthermore, let $\Theta_{\text{sub}}=\{0,\ldots,N_{n}\}\times\{0,\ldots,N_{m}\}\times\{0,\ldots,N_{h}\}$ be the state space of the process $\theta_{t}=\Big{(}\theta^{(n)}(t),\theta^{(m)}(t),\theta^{(h)}(t)\Big{)}$ which records the number of open gates at time t. Note that, $N_{x}-\theta^{(x)}(t)$ gives the number of close gates $x$ at time $t$. The \textit{subunit model} takes the following form\\
\[
(S)\hspace{0.2cm}\left\{
\begin{array}{llll}
C\frac{dV_{t}}{dt}=f^{\text{sub}}(\theta_{t},V_{t},t)\\
(\theta_{t})\\
\end{array}\right.
\]
Where \[f^{\text{sub}}(\theta,V,t)=I(t)-g_{\text{L}}\Big{(}V-V_{\text{L}}\Big{)}-g_{\text{Na}}N_{m}^{-3}\Big{(}\theta^{(m)}\Big{)}^{3}N_{h}^{-1}\theta^{(h)}\Big{(}V-V_{\text{Na}}\Big{)}-g_{\text{K}}N_{n}^{-4}\Big{(}\theta^{(n)}\Big{)}^{4}\Big{(}V-V_{\text{K}}\Big{)}\]
We also define the jump rate of the process by 
\begin{align*}
\lambda^{\text{sub}}(\theta,V)=&\Big{(}\alpha_{m}(V)(N_{m}-\theta^{(m)})+\beta_{m}(V)\theta^{(m)}\Big{)}+\Big{(}\alpha_{h}(V)(N_{h}-\theta^{(h)})+\beta_{h}(V)\theta^{(h)}\Big{)}+\\&\Big{(}\alpha_{n}(V)(N_{n}-\theta^{(n)})+\beta_{n}(V)\theta^{(n)}\Big{)}
\end{align*}
The component $V$ represents the membrane potential and is continuous, thus, the Markovian kernel $Q^{\text{sub}}$ is only concerned by the post-jump location of the jump process $\theta$. We suppose that two channels do not change states simultaneously almost surely. For example, the probability of the event of exactly one gate $n$ opens (conditionally on the last jump time being $T_{k}$) is given by
\[
Q^{\text{sub}}\Big{(}(\theta_{T_{k-1}},V_{T_{k}}),\{\theta_{T_{k-1}}+(1,0,0)\}\Big{)}=\frac{\alpha_{n}(V_{T_{k}})(N_{n}-\theta^{(n)}(T_{k-1}))}{\lambda^{\text{sub}}(\theta_{T_{k-1}},V_{T_{k}})}
\]
To summarize, the \textit{subunit model} can be expressed as a PDMP $x^{\text{sub}}_{t}=(\theta_{t},V_{t},t)\in{\Theta_{\text{sub}}\times\mathbb{R}\times
\mathbb{R}_{+}}$ with vector field $f^{\text{sub}}:\Theta_{\text{sub}}\times\mathbb{R}\times
\mathbb{R}_{+}\rightarrow\mathbb{R}$, jump rate $\lambda^{\text{sub}}:\Theta_{\text{sub}}\times\mathbb{R}\rightarrow\mathbb{R}_{+}$, and a Markovian kernel $Q^{\text{sub}}:\Theta_{\text{sub}}\times\mathbb{R}\times\mathcal{B}(\Theta_{\text{sub}})\rightarrow[0,1]$. 

\subsubsection{The channel model}

In the \textit{channel model}, we form groups of four gates to make channels. 
Unlike the \textit{subunit model}, we define independent MJPs $u_{k}^{(\text{Na})}$ for $k=1,\ldots,N_{\text{Na}}$ (respectively $u_{k}^{(\text{K})}$ for $k=1,\ldots,N_{\text{K}}$) to model the sodium (respectively potassium) channels. These independent MJPs follow kinetic scheme given in appendix B. Note that the conducting state (the state that makes the channel active) of sodium (respectively potassium) channels is $\{m_{3}h_{1}\}$ (respectively $\{n_{4}\}$) which corresponds to three open gates $m$ and one open gate $h$ (respectively four open gates $n$). The conductance of the membrane depends on the empirical measure defined by the proportion of active channels and we define the number of active channels at time $t\geq0$ by
\begin{align*}
\theta^{(m_{3}h_{1})}(t)&=\sum_{k=1}^{N_{\text{Na}}}\textbf{1}_{\textbf{$\{m_{3}h_{1}\}$}}\Big{(}u^{(\text{Na})}_{k}(t)\Big{)},\\
\theta^{(n_{4})}(t)&=\sum_{k=1}^{N_{\text{K}}}\textbf{1}_{\textbf{$\{n_{4}\}$}}\Big{(}u^{(\text{K})}_{k}(t)\Big{)}.
\end{align*}
For $i=0,1,2,3$ and $j=0,1$, let $\theta^{(m_{i}h_{j})}$ be the number of channels in state $\{m_{i}h_{j}\}$ and for $k=0,1,2,3,4$, let $\theta^{(n_{k})}$ be the number of channels in state $\{n_{k}\}$. Let 
\[
\Theta_{\text{chan}}=\{\theta\in{}\{0,\ldots,N_{\text{Na}}\}^{8}\times\{0,\ldots,N_{\text{K}}\}^{5} \hspace{0.1cm}:\hspace{0.1cm} \sum_{i=0}^{3}\sum_{j=0}^{1}\theta^{(m_{i}h_{j})}=N_{\text{Na}},\hspace{0.2cm}\sum_{k=0}^{5}\theta^{(n_{k})}=N_{\text{K}}\}
\]
be the state space of the process $\theta_{t}
=\Big{(}(\theta^{(m_{i}h_{j})}(t))_{i,j},(\theta^{(n_{k})}(t))_{k}\Big{)}$
which records the configuration of the channels at time $t$. The \textit{channel model} takes the following form
\[
(C)\hspace{0.2cm}\left\{
\begin{array}{llll}
C\frac{dV_{t}}{dt}=f^{\text{chan}}(\theta_{t},V_{t},t)\\
(\theta_{t})\\
\end{array}\right.
\]
where
\[
f^{\text{chan}}(\theta,V,t)=I(t)-g_{\text{L}}\Big{(}V-V_{\text{L}}\Big{)}-g_{\text{Na}}N_{\text{Na}}^{-1}\theta^{(m_{3}h_{1})}\Big{(}V-V_{\text{Na}}\Big{)}-g_{\text{K}}N_{\text{K}}^{-1}\theta^{(n_{4})}\Big{(}V-V_{\text{K}}\Big{)}. 
\]
A change in the configuration of the channels (which can be observable or not unlike the \textit{subunit model} in which all changes are observable) happens when a gate opens or closes. We define the application $\eta : \Theta_{\text{chan}}\rightarrow\Theta_{\text{sub}}$ which, given a configuration of channels, returns the configuration of the corresponding gates. We have 
\[
\eta(\theta)=\begin{bmatrix}\theta^{(n_{1})}+2\theta^{(n_{2})}+3\theta^{(n_{3})}+4\theta^{(n_{4})}\\	\theta^{(m_{1}h_{0})}+2\theta^{(m_{2}h_{0})}+3\theta^{(m_{3}h_{0})}+\theta^{(m_{1}h_{1})}+2\theta^{(m_{2}h_{1})}+3\theta^{(m_{3}h_{1})} \\\theta^{(m_{0}h_{1})}+\theta^{(m_{1}h_{1})}+\theta^{(m_{2}h_{1})}+\theta^{(m_{3}h_{1})} \end{bmatrix}
\]
The first component of the vector $\eta(\theta)$ contains $\theta^{n}_{\text{open}}$, the number of open gates $n$, the second $\theta^{m}_{\text{open}}$, the number of open gates $m$ and the third $\theta^{h}_{\text{open}}$, the number of open gates $h$. Thus, for $x=m,h,n$, $\theta^{x}_{\text{close}}(t)=N_{x}-\theta^{x}_{\text{open}}(t)$ gives the number of close gates $x$ at time $t$.
We define the jump rate of the \textit{channel model} by
\[
\lambda^{\text{chan}}(\theta,V)=\lambda^{\text{sub}}(\eta(\theta),V)
\]
Where
\begin{align*}
\lambda^{\text{sub}}(\eta(\theta),V)=&\Big{(}\alpha_{m}(V)(N_{m}-\theta^{m}_{\text{open}})+\beta_{m}(V)\theta^{m}_{\text{open}}\Big{)}+\Big{(}\alpha_{h}(V)(N_{h}-\theta^{h}_{\text{open}})+\beta_{h}(V)\theta^{h}_{\text{open}}\Big{)}+\\&\Big{(}\alpha_{n}(V)(N_{n}-\theta^{n}_{\text{open}})+\beta_{n}(V)\theta^{n}_{\text{open}}\Big{)}
\end{align*}
Since $V$ is continuous, the kernel $Q^{\text{chan}}$ is also only concerned by the post-location of the process $\theta$. Define $Q^{\text{chan}}$ as it is classically do in the literature (\cite{rid} p.53 and \cite{mino} p.587)  is computationally expensive because we have more transitions to deal with than in the \textit{subunit model}. We propose to decompose the kernel $Q^{\text{chan}}$ into a product of two kernels. The decomposition is based on the following observation : it is a change in the configuration of the gates that implies a change in the configuration of the channels. Thus, to determine which transition occurs at time $t$ among the 28 transitions given in appendix B, we first determine which gate opens or closes by using the kernel $Q^{\text{sub}}$ with $\lambda^{\text{sub}}(\eta(.),.)$ and then, depending on which gate changes state, we determine a channel transition by using another kernel. For example, suppose that at time $t$ a gate $m$ opens, thus, the possible channel transitions are : $\{m_{0}h_{0}\rightarrow{m}_{1}h_{0}\}$, $\{m_{1}h_{0}\rightarrow{m}_{2}h_{0}\}$, $\{m_{2}h_{0}\rightarrow{m}_{3}h_{0}\}$, $\{m_{0}h_{1}\rightarrow{m}_{1}h_{1}\}$, $\{m_{1}h_{1}\rightarrow{m}_{2}h_{1}\}$, $\{m_{2}h_{1}\rightarrow{m}_{3}h_{1}\}$ and the next transition is one of those. We define six kernels to take into account all the possibilities.\\
Let $L^{m}_{\text{open}}$, $L^{m}_{\text{close}}$, $L^{h}_{\text{open}}$, $L^{h}_{\text{close}}$, $L^{n}_{\text{open}}$, $L^{n}_{\text{close}}$ be kernels defined on $\Theta_{\text{chan}}\times\mathbb{R}\times\mathcal{B}(\Theta_{\text{chan}})$ with values in $[0,1]$ such that $L^{m}_{\text{open}}$ is the kernel which choose a transition as above, $L^{h}_{\text{open}}$ is a kernel which choose a transition among the following ones $\{m_{0}h_{0}\rightarrow{m}_{0}h_{1}\}$, $\{m_{1}h_{0}\rightarrow{m}_{1}h_{1}\}$, $\{m_{2}h_{0}\rightarrow{m}_{2}h_{1}\}$, $\{m_{3}h_{0}\rightarrow{m}_{3}h_{1}\}$ and so on.
For example, the probability of the event of having the transition \{$m_{0}h_{0}\rightarrow{m}_{1}h_{0}$\} (conditional on the last jump time being $T_{k}$) is given by
\begin{align*}
&Q^{\text{chan}}\Big{(}(\theta_{T_{k-1}},V_{T_{k}}),\{\theta_{T_{k-1}}+(-1,+1,0,\ldots,0)\}\Big{)}=
\\&Q^{\text{sub}}\Big{(}(\eta(\theta_{T_{k-1}}),V_{T_{k}}),\{\eta(\theta_{T_{k-1}})+(0,1,0)\}\Big{)}\times{L}^{m}_{\text{open}}\Big{(}(\theta_{T_{k-1}},V_{T_{k}}),\{\theta_{T_{k-1}}+(-1,+1,0,\ldots,0)\}\Big{)}
\end{align*}
Where
\begin{align*}
&Q^{\text{sub}}\Big{(}(\eta(\theta_{T_{k-1}}),V_{T_{k}}),\{\eta(\theta_{T_{k-1}})+(0,1,0)\}\Big{)}=\frac{\alpha_{m}(V_{T_{k}})\theta^{m}_{\text{close}}(T_{k-1})}{\lambda^{\text{sub}}(\eta(\theta_{T_{k-1}}),V_{T_{k}})}\\
&L^{m}_{\text{open}}\Big{(}(\theta_{T_{k-1}},V_{T_{k}}),\{\theta_{T_{k-1}}+(-1,+1,0,\ldots,0)\}\Big{)}=\frac{3\theta^{(m_{0}h_{0})}(T_{k-1})}{\theta^{m}_{\text{close}}(T_{k-1})}
\end{align*}
Finally, the probability of having the transition \{$m_{0}h_{0}\rightarrow{m}_{1}h_{0}$\} is, as expected, given by the rate of this transition multiplied by the number of channels in the state \{$m_{0}h_{0}$\} divided by the total rate.\\ 
For $x\in{E}$, the support $K^{\text{chan}}_{x}$ of the discrete measure of probability $Q^{\text{chan}}(x,.)$ contains at most 28 elements (depending on the current state $x$), thus, in the worst case we have to do 28 "$if-then$" tests to determine the next transition. With the decomposition of $Q^{\text{chan}}$, we have, in the worst case 12 "$if-then$" tests to do. Indeed, for $x\in{E}$ the support $K^{\text{sub}}_{x}$ of the discrete probability $Q^{\text{sub}}(\eta(x),.)$ contains at most six elements, and the support of the probabilities $L^{m}_{\text{open}}(x,.)$, $L^{m}_{\text{close}}(x,.)$, $L^{h}_{\text{open}}(x,.)$, $L^{h}_{\text{close}}(x,.)$, $L^{n}_{\text{open}}(x,.)$, $L^{n}_{\text{close}}(x,.)$ contains also at most six elements (when we deal with a transition of a gate $m$). Therefore, it is computationally cheaper to decompose the kernel.\\ 
Thus, the \textit{channel model} can be expressed as a PDMP $x^{\text{chan}}_{t}=(\theta_{t},V_{t},t)\in{
\Theta_{\text{chan}}\times\mathbb{R}\times\mathbb{R}_{+}}$ with vector field $f^{\text{chan}}:\Theta_{\text{chan}}\times\mathbb{R}\times\mathbb{R}_{+}\rightarrow\mathbb{R}$, jump rate $\lambda^{\text{chan}}:\Theta_{\text{chan}}\times\mathbb{R}\rightarrow\mathbb{R}_{+}$, and a Markovian kernel $Q^{\text{chan}}:\Theta_{\text{chan}}\times\mathbb{R}\times\mathcal{B}(\Theta_{\text{chan}})\rightarrow[0,1]$.

\subsubsection{Explicit form of the flow of the PDMP between two successive jump times}

In this section, we determine the explicit expression of the flow of both models. For $n\geq0$, $t\geq{T_{n}}$ and $x\in\{\text{sub},\text{chan}\}$, the trajectory of the flow $\phi$ on $[T_{n},+\infty[$ is given by the following ODE 
\[\left\{
\begin{array}{llll}
\frac{d\phi(t-T_{n},x_{T_{n}})}{dt}=f^{x}\Big{(}\theta_{T_{n}},\phi(t-T_{n},x_{T_{n}}),t\Big{)}=-a^{x}_{n}\phi(t-T_{n},x_{T_{n}})+b^{x}_{n}+\frac{1}{C}I(t)\\
\phi(0,x_{T_{n}})=V_{T_{n}}
\end{array}\right.\]
where
\begin{align*}
&a^{\text{sub}}_{n}=\frac{1}{C}\left(g_{\text{L}}+g_{\text{Na}}N_{m}^{-3}\Big{(}\theta^{(m)}(T_{n})\Big{)}^{3}N_{h}^{-1}\theta^{(h)}(T_{n})+g_{\text{K}}N_{n}^{-4}\Big{(}\theta^{(n)}(T_{n})\Big{)}^{4}\right)\\
&b^{\text{sub}}_{n}=\frac{1}{C}\left(g_{\text{L}}V_{\text{L}}+g_{\text{Na}}V_{\text{Na}}N_{m}^{-3}\Big{(}\theta^{(m)}(T_{n})\Big{)}^{3}N_{h}^{-1}\theta^{(h)}(T_{n})+g_{\text{K}}V_{\text{K}}N_{n}^{-4}\Big{(}\theta^{(n)}(T_{n})\Big{)}^{4}\right)\\
&a^{\text{chan}}_{n}=\frac{1}{C}\left(g_{\text{L}}+g_{\text{Na}}N_{\text{Na}}^{-1}\theta^{(m_{3}h_{1})}(T_{n})+g_{\text{K}}N_{\text{K}}^{-1}\theta^{(n_{4})}(T_{n})\right)\hspace{3cm}\\
&b^{\text{chan}}_{n}=\frac{1}{C}\left(g_{\text{L}}V_{\text{L}}+g_{\text{Na}}V_{\text{Na}}N_{\text{Na}}^{-1}\theta^{(m_{3}h_{1})}(T_{n})+g_{\text{K}}V_{\text{K}}N_{\text{K}}^{-1}\theta^{(n_{4})}(T_{n})\right)
\end{align*}
Then, the flow is
\[\phi_{\theta_{T_{n}}}(t-T_{n},V_{T_{n}})=e^{-a^{x}_{n}(t-T_{n})}\Big{[}V_{T_{n}}+\frac{b^{x}_{n}}{a^{x}_{n}}(e^{a^{x}_{n}(t-T_{n})}-1)+\frac{1}{C}\int_{T_{n}}^{t}e^{a^{x}_{n}(s-T_{n})}I(s)ds\Big{]}
\]
For both models we consider that the stimulation $I$  takes the form $I(t)=K\textbf{1}_{[t_{1},t_{2}]}(t)$ with $K>0$ and $t,t_{1},t_{2}\in{\mathbb{R}_{+}}$.

\section{Simulations}

We now proceed to the simulations of the \textit{subunit model} and the \textit{channel model} by using the thinning procedure of section 3.2. Firstly, we explicit the three bounds for both models. Secondly, we numerically compare the efficiency of the bounds in term of reject and simulation time. Finally, we use the thinning procedure to compute a variable of biological interest for both models.

\subsection{Determination of the jump rate bounds}

For simplicity of presentation, we do not distinguish the flows of the \textit{subunit model} and the \textit{channel model}, one has to use $a^{sub}$ and $b^{sub}$ for the \textit{subunit model} and $a^{chan}$ and $b^{chan}$ for the \textit{channel model}. The determination of the bounds relies on the fact that $\alpha_{n}$, $\alpha_{m}$, $\beta_{h}$ are increasing functions, $\beta_{n}$, $\beta_{m}$, $\alpha_{h}$ are decreasing, and that for $n\geq0$ the flow $\phi(.-T_{n},x_{T_{n}})$ is bounded.

\subsubsection{The global bound}

To determine the \textit{global bound} we use a result in \cite{stried} concerning the \textit{channel model} which state that if $V_{0}\in{[V_{-},V_{+}]}$ then $V_{t}\in{[V_{-},V_{+}]}$ $\forall{t}\geq0$ with $V_{-}=\min\{V_{\text{Na}},V_{\text{K}},V_{\text{L}}\}$
and $V_{+}=\max\{V_{\text{Na}},V_{\text{K}},V_{\text{L}}\}$. By using the monotony of the opening and closing rate functions, we find
\[\tilde{\lambda}^{\text{glo}}=N_{m}\alpha_{m}(V_{\text{Na}})+N_{h}\beta_{h}(V_{\text{Na}})+N_{n}\alpha_{n}(V_{\text{Na}})\]
The result in \cite{stried} is also applicable to the \textit{subunit model} and leads to the same expression of the \textit{global bound} for this model.

\subsubsection{The local bound}

Let $n\geq0$ and $t\geq{T_{n}}$. We note $\tilde{\lambda}_{n}=\sup_{t\geq{T_{n}}}\lambda\left(\psi(t-T_{n},x_{T_{n}})\right)$. To determine the constant $\tilde{\lambda}_{n}$ 
we write the flow as follows
\[
\phi(t-T_{n},x_{T_{n}})=f_{n}(t)+g_{n}(t)
\]
where
\begin{align*}
&f_{n}(t)=e^{-a_{n}(t-T_{n})}\Big{(}V_{T_{n}}+\frac{b_{n}}{a_{n}}(e^{a_{n}(t-T_{n})}-1)\Big{)}\\
&g_{n}(t)=e^{-a_{n}(t-T_{n})}\frac{1}{C}\int_{T_{n}}^{t}e^{a_{n}(s-T_{n})}I(s)ds
\end{align*}
The purpose is to determine a lower and an upper bound of the flow. We have $a_{n}>0$, $b_{n}$ may be negative or non-negative, and $f_{n}$ is monotone.
By using the fact that $\forall{t}\geq0$, $I(t)\leq{K}$  we find that
\[ \underline{V}_{T_{n}}\leq{}\phi(t-T_{n},x_{T_{n}})\leq{\overline{V}_{T_{n}}}\tag{2}\]
where $\overline{V}_{T_{n}}=V_{T_{n}}\vee\frac{b_{n}}{a_{n}}+\frac{K}{Ca_{n}}$, and $ 
\underline{V}_{T_{n}}=V_{T_{n}}\wedge\frac{b_{n}}{a_{n}}.$
Then, by using the monotony of the opening and closing rate functions we obtain
\begin{align*}
\tilde{\lambda}_{n}=&\Big{(}\alpha_{m}(\overline{V}_{T_{n}})(N_{m}-\theta^{m}_{\text{open}}(T_{n}))+\beta_{m}(\underline{V}_{T_{n}})\theta^{m}_{\text{open}}(T_{n})\Big{)}+\Big{(}\alpha_{h}(\underline{V}_{T_{n}})(N_{h}-\theta^{h}_{\text{open}}(T_{n}))+\\&\beta_{h}(\overline{V}_{T_{n}})\theta^{h}_{\text{open}}(T_{n})\Big{)}+\Big{(}\alpha_{n}(\overline{V}_{T_{n}})(N_{n}-\theta^{n}_{\text{open}}(T_{n}))+\beta_{n}(\underline{V}_{T_{n}})\theta^{n}_{\text{open}}(T_{n})\Big{)}.
\end{align*}
The expression of the \textit{local bound} is the same for the \textit{channel} and \textit{subunit model} but the Markov chain $\theta$ is different.

\subsubsection{The optimal bound}

In the case of the \textit{optimal bound}, we consider two partitions of $[T_{n},+\infty[$. The first one is the same as in section 4.3 which is noted, for fixed $\epsilon>0$, $(\mathcal{P}_{k}^{T_{n},\epsilon})_{k\in{\mathbb{N}}}$. We recall that, for $k\in{\mathbb{N}}$, $\mathcal{P}_{k}^{T_{n},\epsilon}=[T_{n}+k\epsilon,T_{n}+(k+1)\epsilon[$. 
We precise that the integrated \textit{optimal bound} is given for $t\geq{T_{n}}$ by
\[
\tilde{\Lambda}^{\epsilon}_{n}(t)=
\sum_{k\geq0}\sup_{s\in{\mathcal{P}_{k}^{T_{n},\epsilon}}}
\lambda\Big{(}\psi(s-T_{n},x_{T_{n}})\Big{)}\Big{[}(k+1)\epsilon\wedge(t-T_{n})-k\epsilon\wedge(t-T_{n})\Big{]}
\]
and that the inverse is given by
\[
\left(\tilde{\Lambda}^{\epsilon}_{n}\right)^{-1}(t)=\sum_{p\geq0}\Big{(}\frac{t-\epsilon\sum_{k=0}^{p-1}\sup_{s\in{\mathcal{P}_{k}^{T_{n},\epsilon}}}
\lambda\Big{(}\psi(s-T_{n},x_{T_{n}})\Big{)}}{\sup_{s\in{\mathcal{P}_{p}^{T_{n},\epsilon}}}
\lambda\Big{(}\psi(s-T_{n},x_{T_{n}})\Big{)}}+T_{n}+p\epsilon\Big{)}
\textbf{1}_{[\kappa_{p-1},\kappa_{p}[}(t)
\]
where 
$\kappa_{p}=\epsilon\sum_{k=0}^{p}\sup_{s\in{\mathcal{P}_{k}^{T_{n},\epsilon}}}\lambda\Big{(}\psi(s-T_{n},x_{T_{n}})\Big{)}$
and, by convention, $\kappa_{-1}=0$. For $k\in{\mathbb{N}}$, we have
\begin{align*}
\sup_{s\in{\mathcal{P}_{k}^{T_{n},\epsilon}}}\lambda\Big{(}\psi(s-T_{n},x_{T_{n}})\Big{)}=&\Big{(}\alpha_{m}(\overline{V}^{k,\epsilon}_{T_{n}})(N_{m}-\theta^{m}_{open}(T_{n}))+\beta_{m}(\underline{V}^{k,\epsilon}_{T_{n}})\theta^{m}_{open}(T_{n})\Big{)}+\\&\Big{(}\alpha_{h}(\underline{V}^{k,\epsilon}_{T_{n}})(N_{h}-\theta^{h}_{open}(T_{n}))+\beta_{h}(\overline{V}^{k,\epsilon}_{T_{n}})\theta^{h}_{open}(T_{n})\Big{)}+\\&\Big{(}\alpha_{n}(\overline{V}^{k,\epsilon}_{T_{n}})(N_{n}-\theta^{n}_{open}(T_{n}))+\beta_{n}(\underline{V}^{k,\epsilon}_{T_{n}})\theta^{n}_{open}(T_{n})\Big{)}
\end{align*}
Where
\begin{align*} 
&\overline{V}^{k,\epsilon}_{T_{n}}=f_{n}\left(T_{n}+k\epsilon\right)\vee{}f_{n}\left(T_{n}+(k+1)\epsilon\right)
+e^{-a_{n}k\epsilon}\int_{T_{n}}^{T_{n}+(k+1)\epsilon}e^{a_{n}(s-T_{n})}I(s)ds,\\
&\underline{V}^{k,\epsilon}_{T_{n}}=f_{n}\left(T_{n}+k\epsilon\right)\wedge{}f_{n}\left(T_{n}+(k+1)\epsilon\right)
+e^{-a_{n}(k+1)\epsilon}\int_{T_{n}}^{T_{n}+k\epsilon}e^{a_{n}(s-T_{n})}I(s)ds.
\end{align*}
The second partition is obtained for $P=\{0,1\}$ and is noted $(\mathcal{Q}_{k}^{T_{n},\epsilon})_{k\in{P}}$ where $\mathcal{Q}_{0}^{T_{n},\epsilon}=[T_{n},T_{n}+\epsilon[$ and $\mathcal{Q}_{1}^{T_{n},\epsilon}=[T_{n}+\epsilon,+\infty[$.
The integrated \textit{optimal bound} is
\[
\tilde{\Lambda}^{\epsilon}_{n}(t)
=\overline{\lambda}^{\epsilon}_{n}\Big{[}\epsilon\wedge(t-T_{n})\Big{]}+\tilde{\lambda}_{n}\Big{[}(t-T_{n})-\epsilon\wedge(t-T_{n})\Big{]}
\]
Where $\tilde{\lambda}_{n}$ is the previous \textit{local bound} and
\[
\overline{\lambda}^{\epsilon}_{n}=\sup_{t\in{\mathcal{Q}_{0}^{T_{n},\epsilon}}}\lambda\Big{(}\psi(t-T_{n},x_{T_{n}})\Big{)}.
\]
The inverse is given by
\[
\Big{(}\tilde{\Lambda}^{\epsilon}_{n}\Big{)^{-1}}(s)=\Big{(}\frac{s}{\overline{\lambda}^{\epsilon}_{n}}+T_{n}\Big{)}\textbf{1}_{[0,
\epsilon\overline{\lambda}^{\epsilon}_{n}[}(s)+\Big{(}\frac{s-\epsilon\overline{\lambda}^{\epsilon}_{n}}{\tilde{\lambda}_{n}}+T_{n}+\epsilon\Big{)}\textbf{1}_{[\epsilon\overline{\lambda}^{\epsilon}_{n},+\infty[}(s)
\]
Once again, the expression of the \textit{optimal bound} is the same for both models but the Markov chain is different. We precise that we used the \textit{local bound} to define the \textit{optimal bound} with the partition $(\mathcal{Q}_{k}^{T_{n},\epsilon})_{k\in{\{0,1\}}}$. \\
Note that, for $n\geq1$, it is possible to define an $\epsilon_{n}$ which is "adapted" to the inter jump time $T_{n+1}-T_{n}$.
To determine such an $\epsilon_{n}$, we use the bounds of the flow in inequality $(2)$ to define a lower \textit{local bound} $\underline{\lambda}_{n}$ such that
\begin{align*}
\underline{\lambda}_{n}=&\Big{(}\alpha_{m}(\underline{V}_{T_{n}})(N_{m}-\theta^{m}_{\text{open}}(T_{n}))+\beta_{m}(\overline{V}_{T_{n}})\theta^{m}_{\text{open}}(T_{n})\Big{)}+\Big{(}\alpha_{h}(\overline{V}_{T_{n}})(N_{h}-\theta^{h}_{\text{open}}(T_{n}))+\\&\beta_{h}(\underline{V}_{T_{n}})\theta^{h}_{\text{open}}(T_{n})\Big{)}+\Big{(}\alpha_{n}(\underline{V}_{T_{n}})(N_{n}-\theta^{n}_{\text{open}}(T_{n}))+\beta_{n}(\overline{V}_{T_{n}})\theta^{n}_{\text{open}}(T_{n})\Big{)}.
\end{align*}
We note $(\underline{T}_{n})_{n\geq0}$ the corresponding point process, then we have
\[
0.05=\mathbb{P}(T_{n+1}-T_{n}>\epsilon_{n})\leq{\mathbb{P}(\underline{T}_{n+1}-\underline{T}_{n}>\epsilon_{n})}=e^{-\epsilon_{n}\underline{\lambda}_{n}}
\]
and we take 
\[
\epsilon_{n}=\frac{-\log(0.05)}{\underline{\lambda}_{n}}.
\]
Note that the value $\epsilon_{n}$ is in fact adapted to the inter jump time $\underline{T}_{n+1}-\underline{T}_{n}$.

\subsection{Numerical results}

In this section, we numerically compare the three different jump rate bounds and we use thinning algorithm to simulate a variable of biological interest, the spiking times.
  
\subsubsection{Numerical comparison of the jump rate bounds}

In this part, we first show trajectories of the two stochastic Hodgkin-Huxley models obtained with the thinning method using the \textit{optimal bound}. Then, we collect in several tables and graphs the results concerning the simulation time and the rate of acceptance of both models.

In the sequel, for $\epsilon>0$, the \textit{optimal}-$\mathcal{Q}^{\epsilon}$ (respectively \textit{optimal}-$\mathcal{P}^{\epsilon}$) bound denotes the optimal bound using the partition $(\mathcal{Q}_{k}^{T_{n},\epsilon})_{k\in{\{0,1\}}}$ (respectively $(\mathcal{P}_{k}^{T_{n},\epsilon})_{k\in{\mathbb{N}}}$). All numerical values are obtained from a classical Monte Carlo method with 100 000 trials. Parameters of the models are given in Appendix \ref{parameters} - \ref{Na_K}. We denote by $N_{\text{chan}}$ the common number of sodium and potassium channels, $N_{\text{chan}}=N_{\text{Na}}=N_{\text{K}}$. The input current is $I(t)=30\textbf{1}_{[1,2]}(t)$. The simulation time represents the time needed to simulate one path of the PDMP on [0,10]. 

\begin{figure}
  \begin{subfigure}[b]{.5\linewidth}
    \includegraphics[width=.9\textwidth]{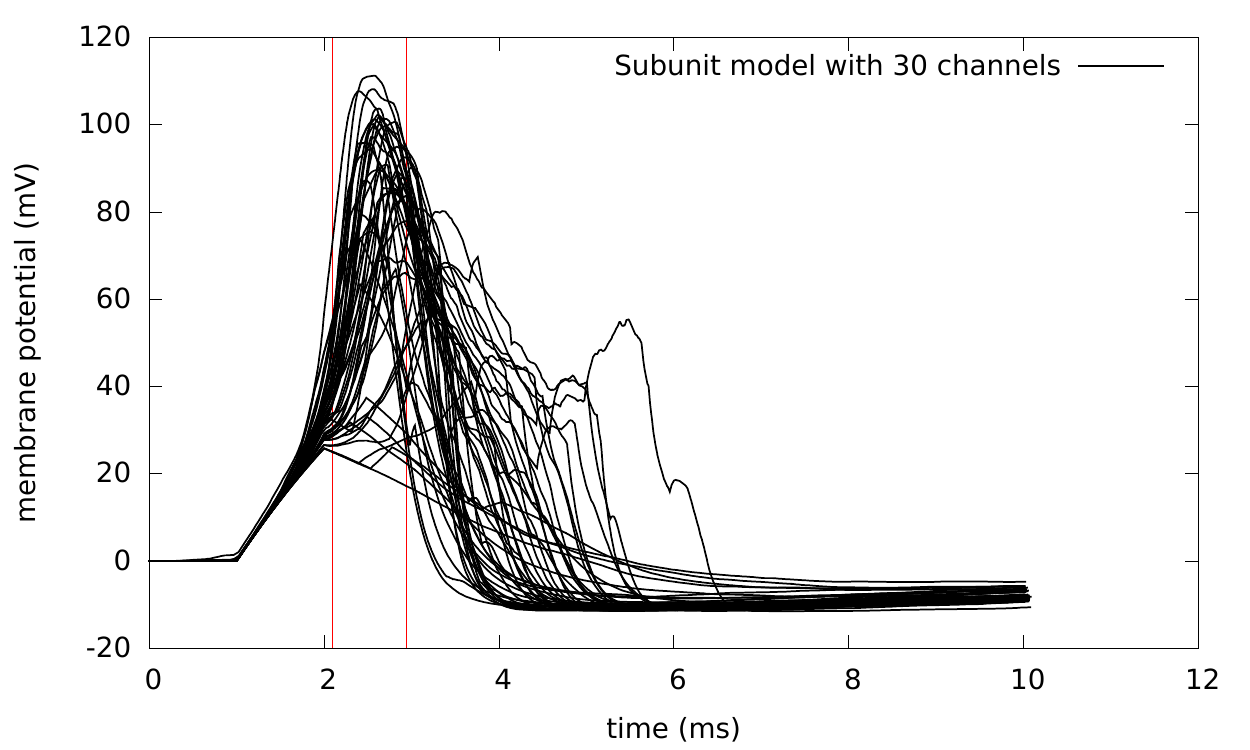}
  \end{subfigure} 
  \begin{subfigure}[b]{.5\linewidth}
    \includegraphics[width=.9\textwidth]{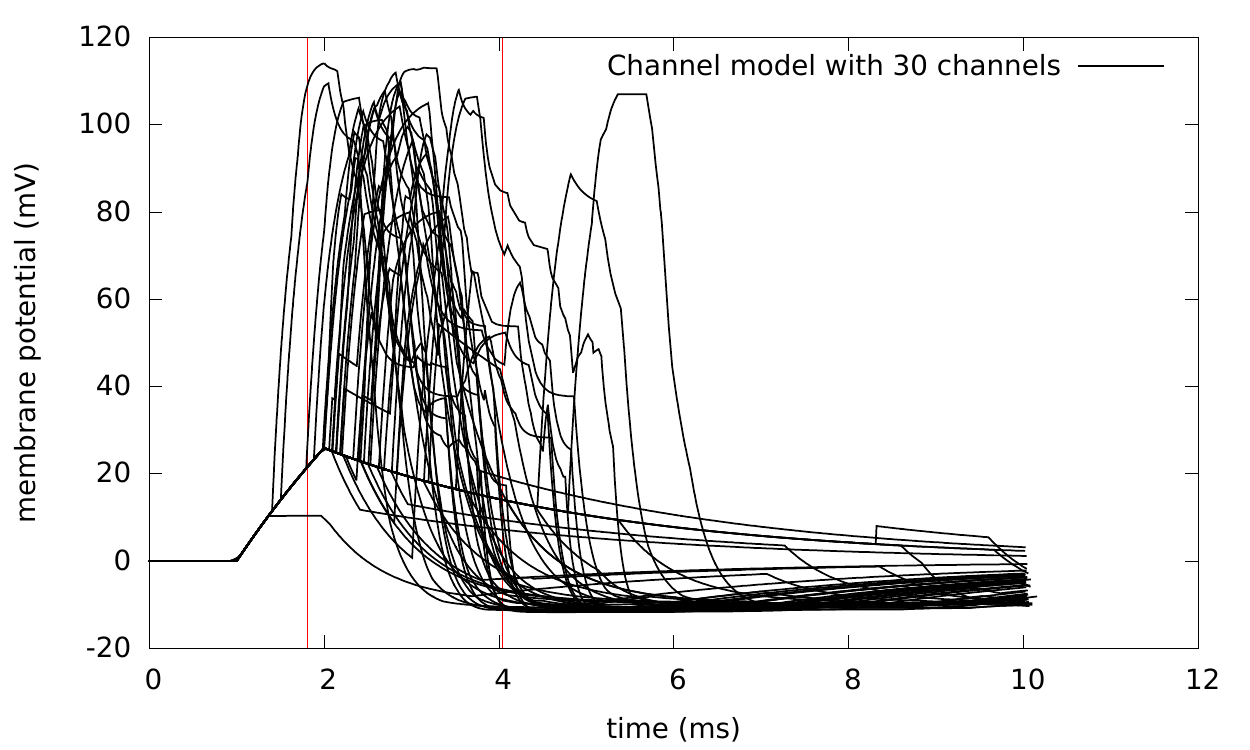}
  \end{subfigure}
  \vspace{0.4cm}
  \begin{subfigure}[b]{.5\linewidth}
    \includegraphics[width=.9\textwidth]{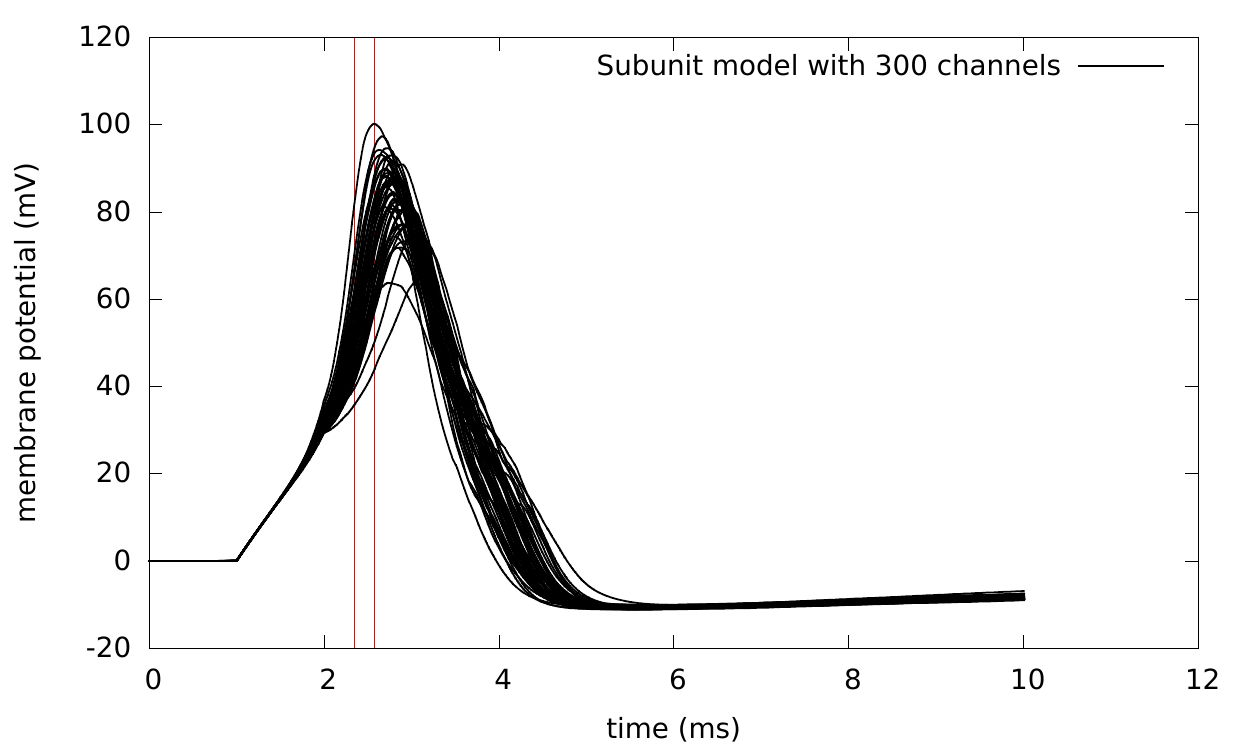}
  \end{subfigure}
  \begin{subfigure}[b]{.5\linewidth}
    \includegraphics[width=.9\textwidth]{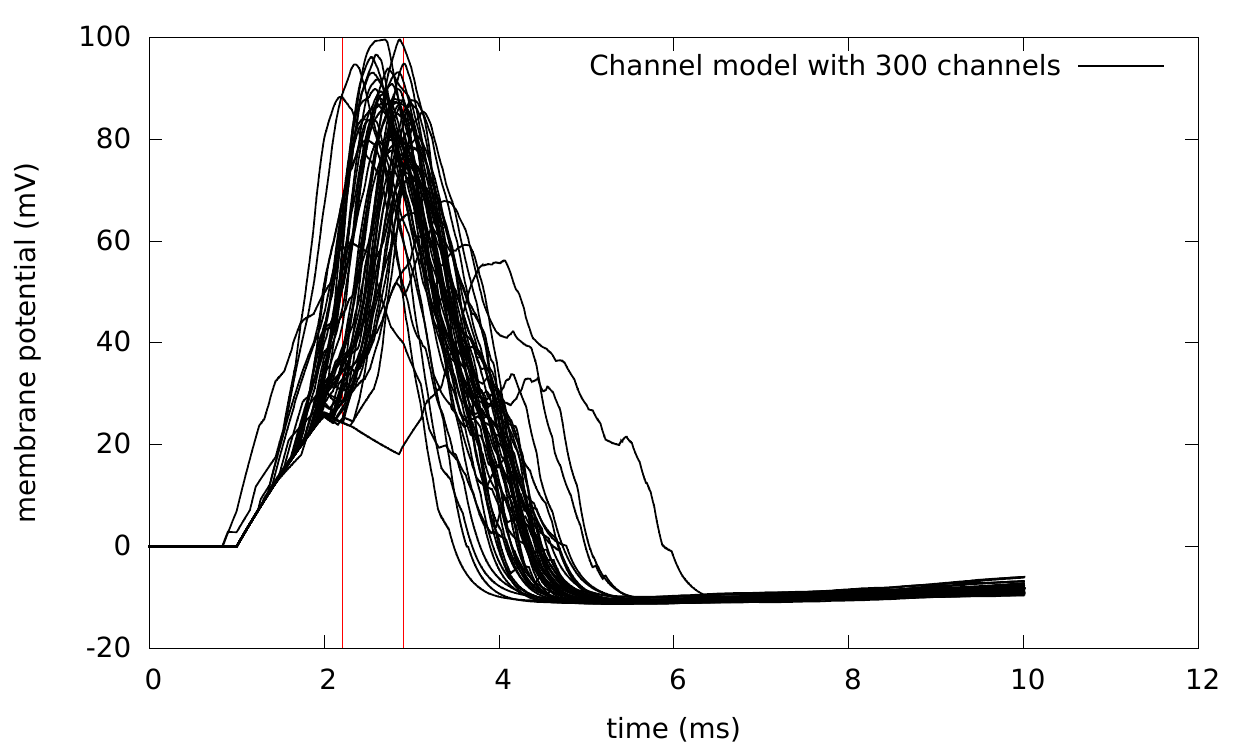}
  \end{subfigure}
  \vspace{0.4cm}
  \begin{subfigure}[b]{.5\linewidth}
    \includegraphics[width=.9\textwidth]{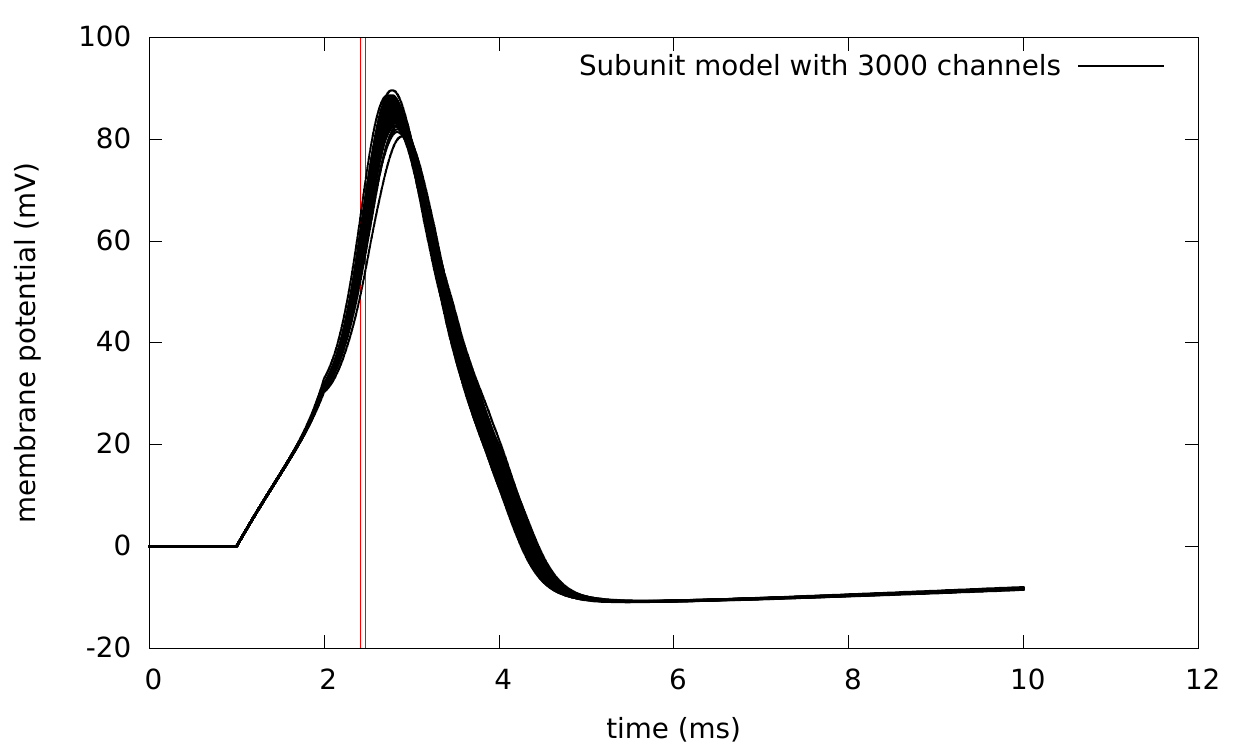}
  \end{subfigure} 
  \begin{subfigure}[b]{.5\linewidth}
    \includegraphics[width=.9\textwidth]{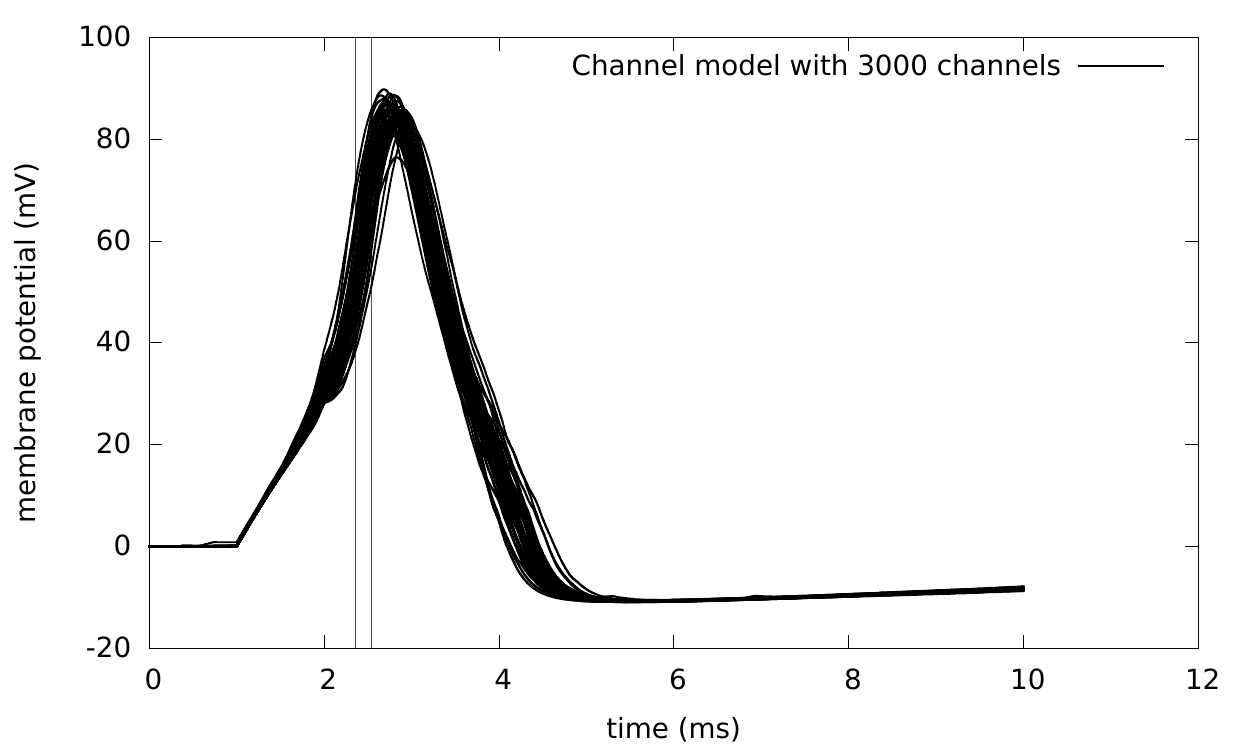}
  \end{subfigure}
  \caption{First column : \textit{subunit model}. Second column : \textit{channel model}. Vertical lines are the standard deviation of the spiking times (see section 6.2.2).}
\label{fig:2}
\end{figure}

Each lines of Figure 2 shows fifty trajectories of the \textit{subunit} and the \textit{channel model} with a different number of channels, $N_{\text{chan}}=30,300,3000$. It allows to see the different behaviours of the two models. 
In each lines, we see that the behaviour of the \textit{channel model} is more erratic than the \textit{subunit model} one (except for the third line where the two models have approximately the same behaviour). Differences in trajectories are mainly explained by two distinct modelling approaches of the conductance of the membrane. In the \textit{subunit model}, we consider that the conductance at time $t$ depends on the fraction of open gates at time $t$, thus, the equation of the voltage changes rapidly at the same time as the state of the gates. 
In the \textit{channel model}, the conductance at time $t$ depends on the fraction of active channels at time $t$, therefore, a change in the state of the gates may not implies a change in the voltage's equation. Thus, the dynamic of the membrane potential change less than in the first case and trajectories are more irregular. We also see that more the number of channels is high more the differences in trajectories are small. It illustrates a result in \cite{wainfluid} where the authors showed that the deterministic limit when the number of channels goes to infinity of both models are equivalent. However, it seems that the convergence speed is not the same.

\begin{figure}
  \begin{subfigure}[b]{.5\linewidth}
    \includegraphics[width=\textwidth]{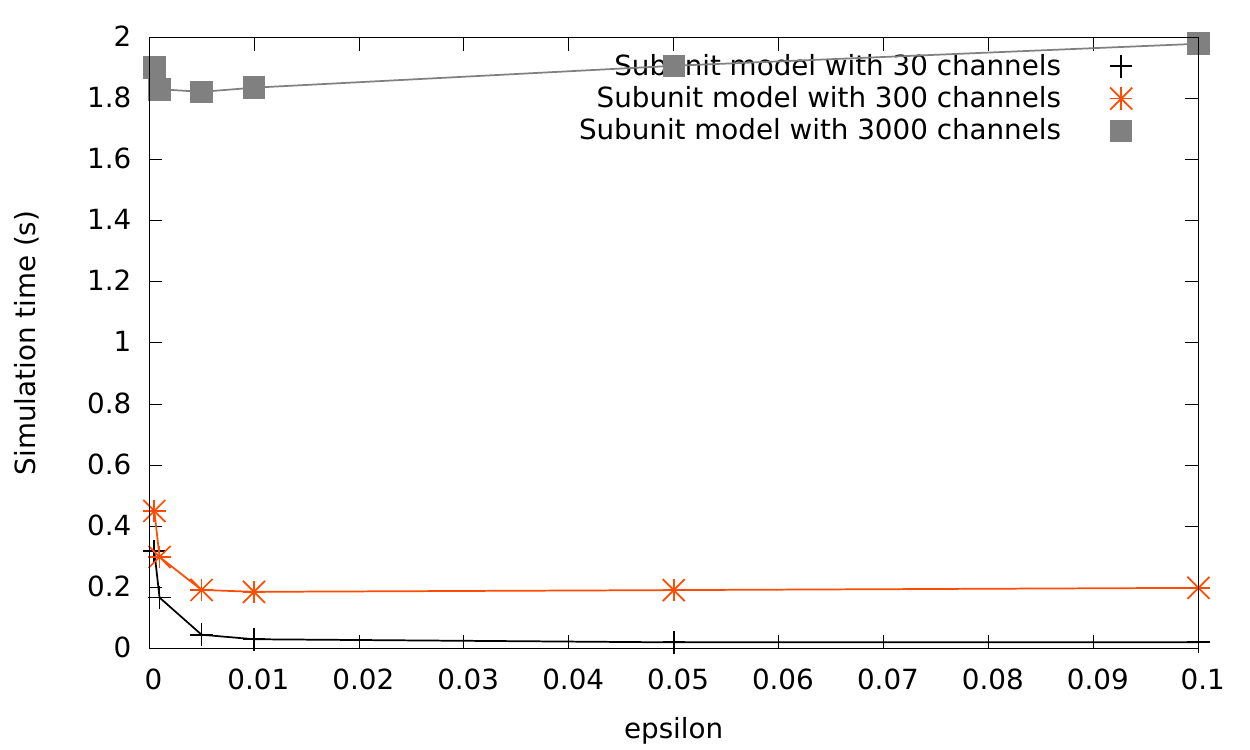}
  \end{subfigure} 
  \begin{subfigure}[b]{.5\linewidth}
    \includegraphics[width=\textwidth]{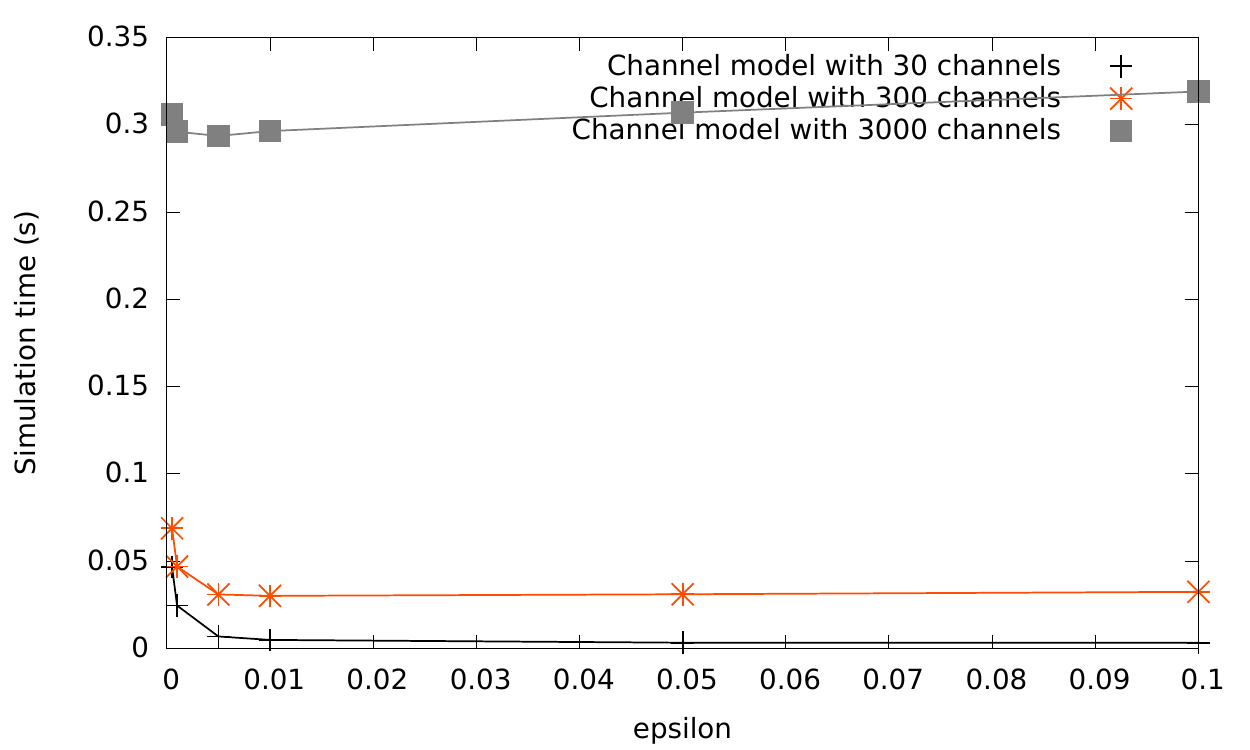}
  \end{subfigure}
  \vspace{0.4cm}
  \begin{subfigure}[b]{.5\linewidth}
    \includegraphics[width=\textwidth]{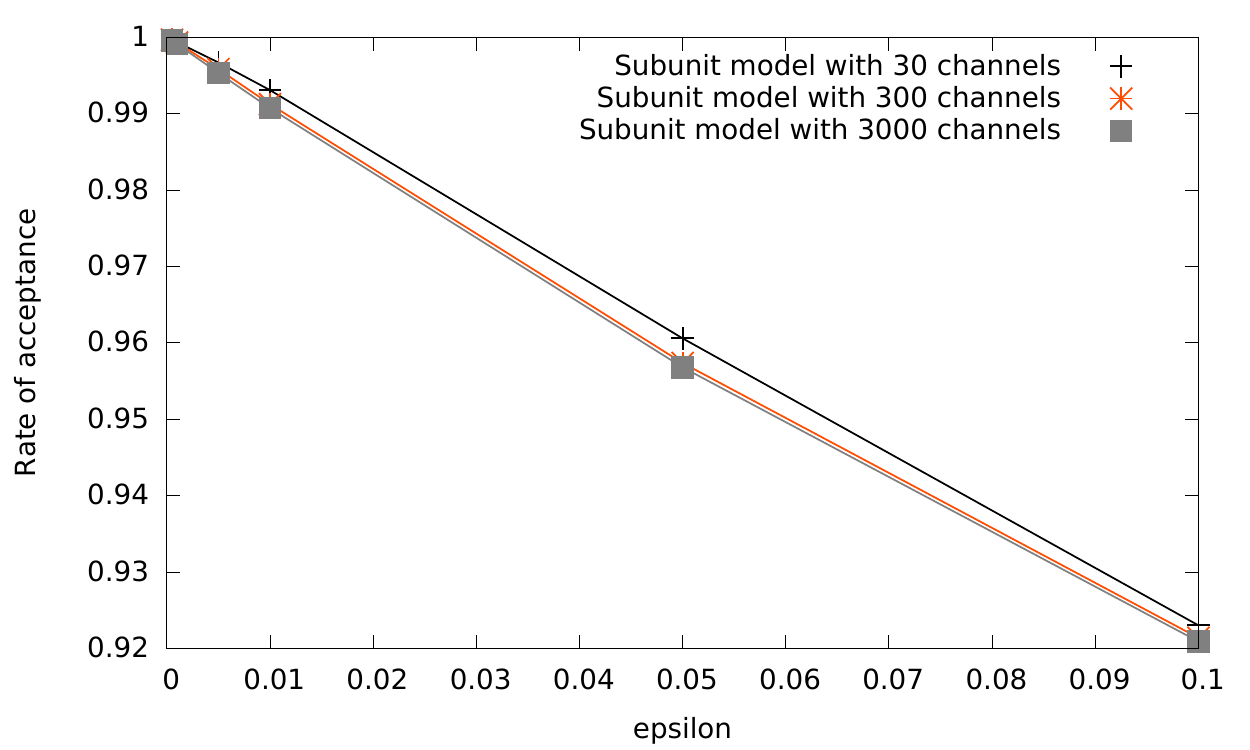}
  \end{subfigure}
  \begin{subfigure}[b]{.5\linewidth}
    \includegraphics[width=\textwidth]{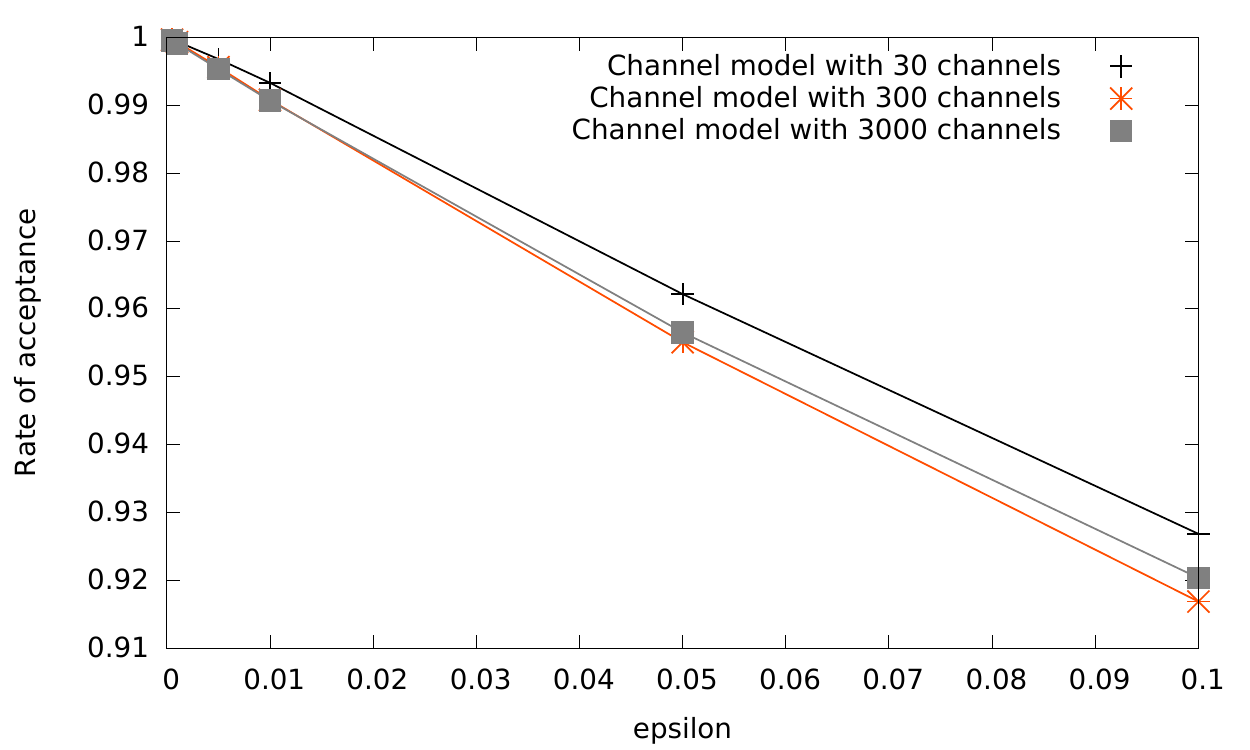}
  \end{subfigure}
  \caption{Simulation time and rate of acceptance with the optimal-$\mathcal{P}^{\epsilon}$ bound as a function of the parameter $\epsilon$.}
\label{fig:3}
\end{figure}

\begin{figure}
  \begin{subfigure}[b]{.5\linewidth}
    \includegraphics[width=\textwidth]{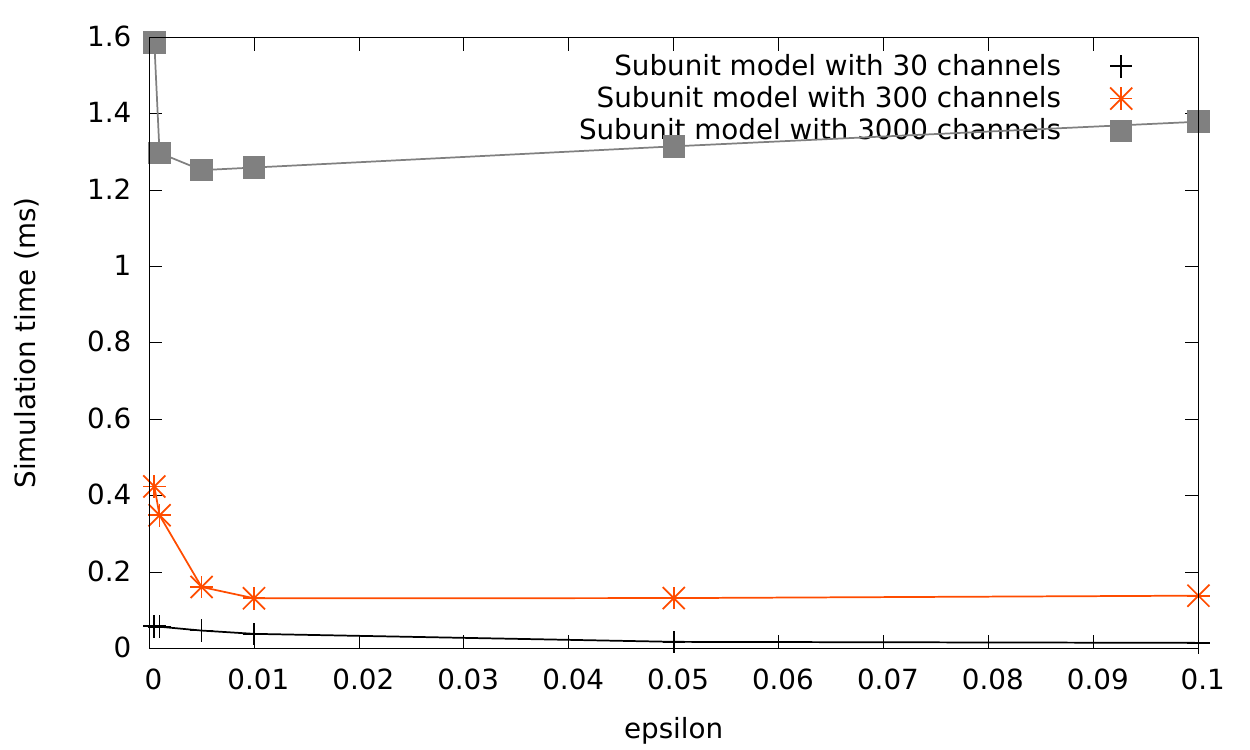}
  \end{subfigure} 
  \begin{subfigure}[b]{.5\linewidth}
    \includegraphics[width=\textwidth]{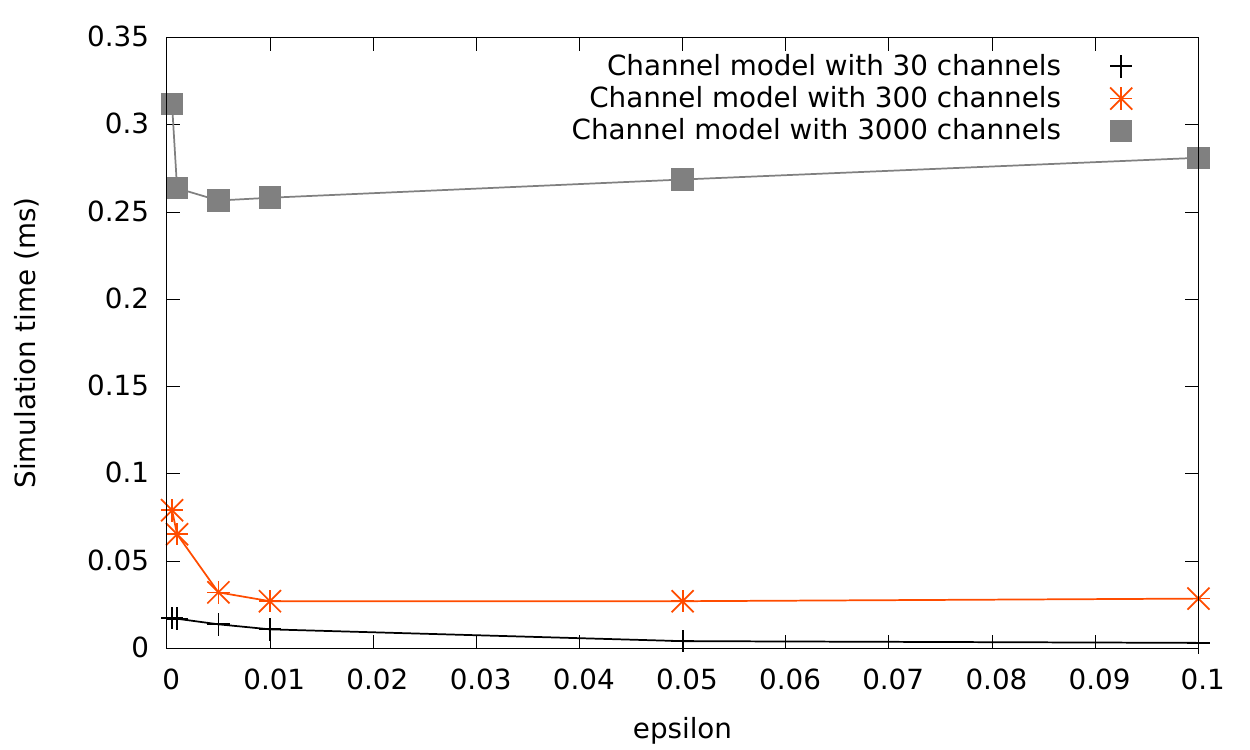}
  \end{subfigure}
  \vspace{0.4cm}
  \begin{subfigure}[b]{.5\linewidth}
    \includegraphics[width=\textwidth]{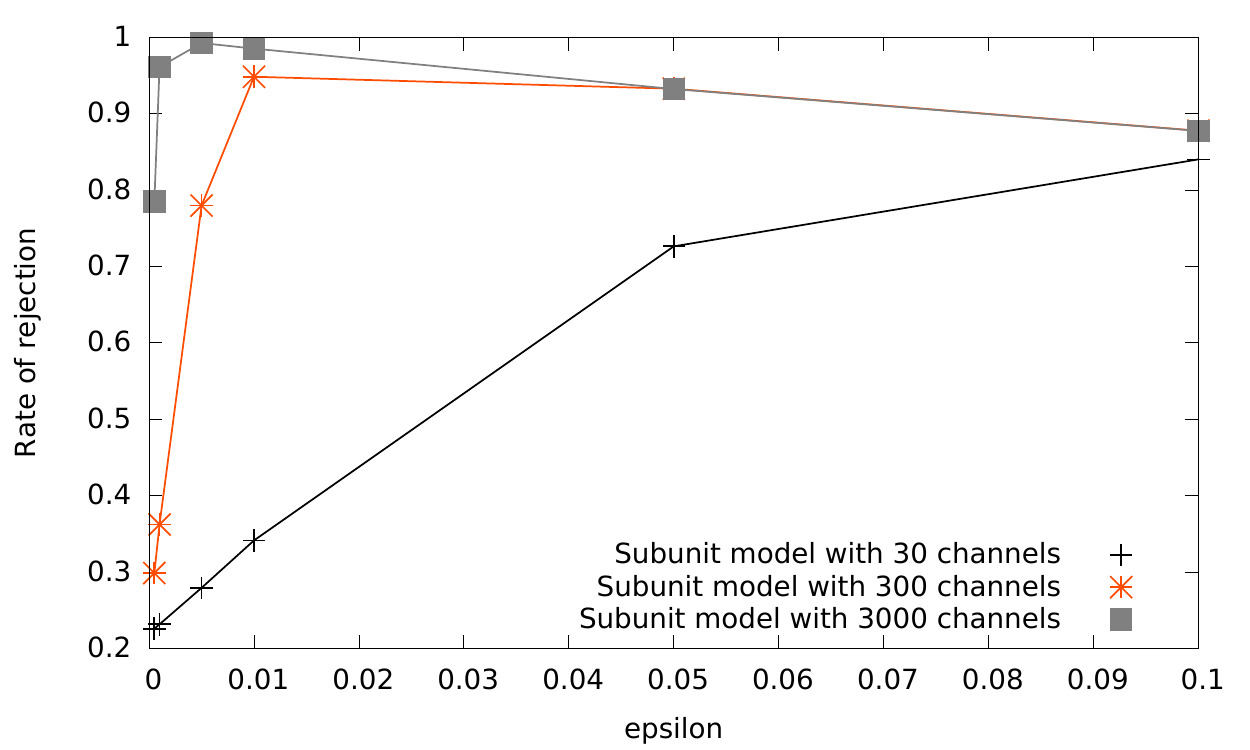}
  \end{subfigure}
  \begin{subfigure}[b]{.5\linewidth}
    \includegraphics[width=\textwidth]{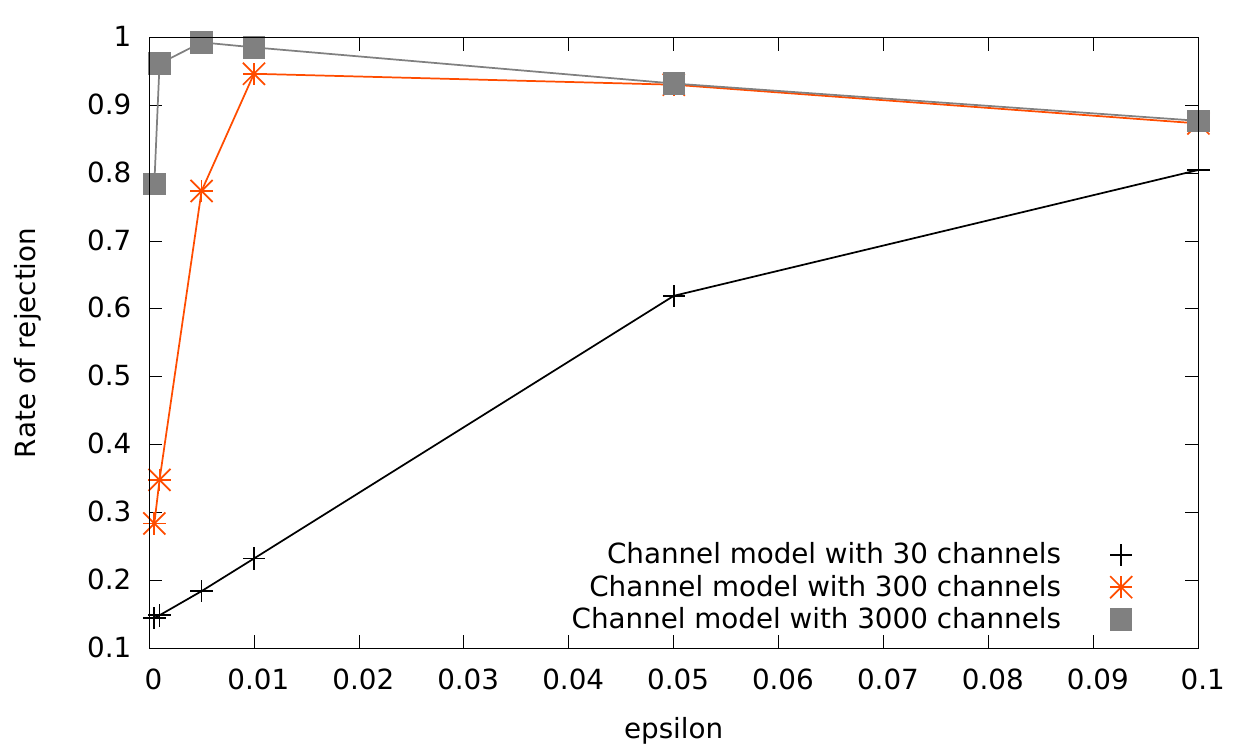}
  \end{subfigure}
  \caption{Simulation time and rate of acceptance with the optimal-$\mathcal{Q}^{\epsilon}$ bound as a function of the parameter $\epsilon$.}
\label{fig:4}
\end{figure}

Concerning the \textit{optimal}-$\mathcal{P}^{\epsilon}$ bound, we see on Figure \ref{fig:3} that in both models, more $\epsilon$ is small less rejected points are. It illustrates the fact that $\tilde{N}^{\epsilon}$ converge to $N$ when $\epsilon$ goes to 0 (proposition 4.4). Figure 3 also shows that, for fixed $N_{\text{chan}}$, the simulation time varies with $\epsilon$. For both models, the value of $\epsilon$ which minimize the simulation time is inversely proportional to the parameter $N_{\text{chan}}$. Let $\epsilon(N_{\text{chan}})$ be that optimal value of $\epsilon$. For increasing $\epsilon>\epsilon(N_{\text{chan}})$, the rate of acceptance decreases, thus, we have to simulate more and more uniform pseudo-random variables and the simulation time increases. For decreasing $\epsilon<\epsilon(N_{\text{chan}})$, the rate of acceptance increases but the simulation time too because of the increasing number of iterations needed to compute the integrated jump rate bound and its inverse. Thus, one has to take a small (respectively large) $\epsilon$ when the jumps frequency is high (respectively low).\\
We see on Figure 4 that more $\epsilon$ is small more the rate of acceptance of the \textit{optimal}-$\mathcal{Q}^{\epsilon}$ bound is close to the one of the \textit{local bound}. Note that the value of $\epsilon$ which maximise the rate of acceptance is the same which minimize the simulation time. As in the case of the \textit{optimal}-$\mathcal{P}^{\epsilon}$ bound, the optimal value of $\epsilon$ is inversely proportional to $N_{\text{chan}}$. For decreasing $\epsilon<\epsilon(N_{\text{chan}})$, the rate of acceptance decreases and the simulation time increases because we use the \textit{local bound} $\tilde{\lambda}_{n}$ instead of the smaller bound $\overline{\lambda}_{n}^{\epsilon}$ (see section 6.1.3). For increasing $\epsilon>\epsilon(N_{\text{chan}})$, the rate of acceptance decreases and the simulation time increases because the bound $\overline{\lambda}_{n}^{\epsilon}$ becomes bigger and bigger. \\
By comparing  the \textit{optimal}-$\mathcal{Q}^{\epsilon(N_{\text{chan}})}$ and the \textit{optimal}-$\mathcal{P}^{\epsilon(N_{\text{chan}})}$ bound we see that the first one is the most efficient in term of simulation time, it is also the simplest to implement. However, this bound does not exist when the jump rate or the flow is not bounded. In this case, one may use the \textit{optimal}-$\mathcal{P}^{\epsilon(N_{\text{chan}})}$ bound which is efficient too but more complex to implement.\\
From figure 3 and 4, we see that for both the \textit{optimal}-$\mathcal{Q}^{\epsilon}$ and the \textit{optimal}-$\mathcal{P}^{\epsilon}$ bound the best simulation time is achieved for $\epsilon(30)=0.1$, $\epsilon(300)=0.01$ and $\epsilon(3000)=0.005$. We saw in sections 5.2.1 and 5.2.2 that the \textit{subunit model} and the \textit{channel model} share the same jump rate. For both models, the maximum value of the inter-jump times is of order $10^{-1}$ for $N_{\text{chan}}=30$, $10^{-2}$ for $N_{\text{chan}}=300$ and $10^{-3}$ for $N_{\text{chan}}=3000$. It coincides with the values $\epsilon(N_{\text{chan}})$ which, in this case, confirm that the optimal simulation time is obtained for $\epsilon$ of order $\max_{n}|T_{n+1}-T_{n}|$. 

\begin{table}
\caption{Simulation time and rate of acceptance for $N_{\text{chan}}=30$. The lines ODE represent the algorithm in \cite{rid} with $h=10^{-3}$ for both \textit{subunit model} and \textit{channel model}}
\label{tab:1}      
\setlength{\extrarowheight}{4 pt}
\begin{tabular}{|c|c|c|c|}
\hline
Model & Bound & simulation time (sec) & rate of acceptance\\
\hline
 & Optimal-$\mathcal{Q}^{\epsilon_{n}}$ & 0,003$\hspace{0.4cm}(\pm8.10^{-7})$ & 0,857$\hspace{0.4cm}(\pm2.10^{-3})$ \\
 \cline{2-4}
Channel & Local & 0,008$\hspace{0.4cm}(\pm6.10^{-6})$ & 0,141$\hspace{0.4cm}(\pm2.10^{-3})$  \\
\cline{2-4}
 & Global & 0,012$\hspace{0.4cm}(\pm3.10^{-6})$ & 0,065$\hspace{0.4cm}(\pm6.10^{-5})$ \\
  \cline{2-4}
 & ODE & 0.009$\hspace{0.4cm}(\pm1.10^{-7})$ &  \\
\hline\hline
 & Optimal-$\mathcal{Q}^{\epsilon_{n}}$ & 0,016$\hspace{0.4cm}(\pm1.10^{-6})$ & 0,88$\hspace{0.4cm}(\pm1.10^{-3})$ \\
 \cline{2-4}
Subunit & Local & 0,050$\hspace{0.4cm}(\pm2.10^{-4})$ & 0,22$\hspace{0.4cm}(\pm1.10^{-3})$  \\
\cline{2-4}
& Global & 0,12$\hspace{0.4cm}(\pm3.10^{-4})$ & 0,061$\hspace{0.4cm}(\pm2.10^{-5})$ \\
  \cline{2-4}
 & ODE & 0.016$\hspace{0.4cm}(\pm2.10^{-7})$ &  \\
\hline
\end{tabular}
\end{table}

\begin{table}
\caption{Simulation time and rate of acceptance for $N_{\text{chan}}=300$. The lines ODE represent the algorithm in \cite{rid} with $h=10^{-4}$ for both \textit{subunit model} and \textit{channel model}.}
\label{tab:2}
\setlength{\extrarowheight}{4 pt}
\begin{tabular}{|c|c|c|c|c|}
\hline
Model & Bound & simulation time (sec) & rate of acceptance\\
\hline
 & Optimal-$\mathcal{Q}^{\epsilon_{n}}$ & 0,030$\hspace{0.4cm}(\pm3.10^{-5})$ & 0,962$\hspace{0.4cm}(\pm9.10^{-5})$ \\
\cline{2-4}
Channel & Local & 0,050$\hspace{0.4cm}(\pm1.10^{-4})$  & 0,223$\hspace{0.4cm}(\pm3.10^{-4})$  \\
\cline{2-4}
 & Global & 0,120$\hspace{0.4cm}(\pm3.10^{-4})$ & 0,062$\hspace{0.4cm}(\pm7.10^{-5})$ \\
 \cline{2-4}
 & ODE & 0.094$\hspace{0.4cm}(\pm1.10^{-5})$ &  \\
\hline\hline
 & Optimal-$\mathcal{Q}^{\epsilon_{n}}$ & 0,148$\hspace{0.4cm}(\pm5.10^{-4})$ & 0,957$\hspace{0.4cm}(\pm9.10^{-5})$ \\
 \cline{2-4}
Subunit & Local & 0,244$\hspace{0.4cm}(\pm1.10^{-3})$ & 0,237$\hspace{0.4cm}(\pm8.10^{-5})$  \\
\cline{2-4}
& Global & 0,322$\hspace{0.4cm}(\pm2.10^{-3})$ & 0,061$\hspace{0.4cm}(\pm1.10^{-5})$ \\
\cline{2-4}
 & ODE & 0.157$\hspace{0.4cm}(\pm1.10^{-5})$ &  \\
\hline
\end{tabular}
\end{table}

\begin{table}
\caption{Simulation time and rate of acceptance for $N_{\text{chan}}=3000$. The lines ODE represent the algorithm in \cite{rid} with $h=10^{-5}$ for both \textit{subunit model} and \textit{channel model}.}
\label{tab:3}
\setlength{\extrarowheight}{4 pt}
\begin{tabular}{|c|c|c|c|}
\hline
Model & Bound & simulation time (sec) & rate of acceptance\\
\hline
 &  optimal-$\mathcal{Q}^{\epsilon_{n}}$ & 0,296$\hspace{0.4cm}(\pm3.10^{-3})$ & 0,965$\hspace{0.4cm}(\pm2.10^{-5})$ \\
\cline{2-4}
Channel & Local & 0,474$\hspace{0.4cm}(\pm6.10^{-3})$  & 0,236$\hspace{0.4cm}(\pm3.10^{-5})$  \\
\cline{2-4}
 & Global & 1,184$\hspace{0.4cm}(\pm2.10^{-2})$ & 0,060$\hspace{0.4cm}(\pm3.10^{-7})$ \\
 \cline{2-4}
 & ODE & 0.940$\hspace{0.4cm}(\pm5.10^{-4})$ &  \\
\hline\hline
 & Optimal-$\mathcal{Q}^{\epsilon_{n}}$  & 1,471$\hspace{0.4cm}(\pm3.10^{-2})$ & 0,964$\hspace{0.4cm}(\pm9.10^{-6})$ \\
\cline{2-4}
Subunit & Local & 2,478$\hspace{0.4cm}(\pm4.10^{-2})$ & 0,238$\hspace{0.4cm}(\pm7.10^{-6})$  \\
\cline{2-4}
 & Global & 3,315$\hspace{0.4cm}(\pm3.10^{-1})$ & 0,060$\hspace{0.4cm}(\pm9.10^{-8})$ \\
\cline{2-4}
 & ODE & 1.567$\hspace{0.4cm}(\pm1.10^{-3})$ &  \\
\hline
\end{tabular}
\end{table}

Tables \ref{tab:1}-\ref{tab:3} show results of numerical computations of the simulation time and of the rate of acceptance of the thinning method for the \textit{global}, \textit{local} and \textit{optimal}-$\mathcal{Q}^{\epsilon_{n}}$ bound using both the \textit{channel} and the \textit{subunit model} with different values of the parameter $N_{\text{chan}}$. For both models and for all the studied values of $N_{\text{chan}}$, the simulation time using the \textit{optimal} bounds ($\mathcal{Q}^{\epsilon(N_{\text{chan}})}$,$\mathcal{P}^{\epsilon(N_{\text{chan}})}$ and $\mathcal{Q}^{\epsilon_{n}}$) is better than the one obtained with both the \textit{global} and \textit{local bound}. Note that the \textit{optimal}-$\mathcal{Q}^{\epsilon_{n}}$ bound is more efficient than the \textit{optimal}-$\mathcal{P}^{\epsilon(N_{\text{chan}})}$ bound to simulate the \textit{subunit model}. Since the computation of $\epsilon_{n}$ requires the computation of the jump rate bound at each iterations, the \textit{optimal}-$\mathcal{Q}^{\epsilon_{n}}$ bound will be more efficient when the jumps frequency is low. Thus, for all studied values of $N_{\text{chan}}$, the \textit{optimal}-$\mathcal{Q}^{\epsilon(N_{\text{chan}})}$ is the most efficient.
\\
The differences of simulation time between the \textit{subunit} and the \textit{channel model} are explained by the fact that the numerical computation of the flow of the \textit{channel model} is cheaper than the one of the \textit{subunit model}. Note that the simulation time using the three bounds (\textit{global}, \textit{local}, \textit{optimal}) increases linearly as a function of $N_{\text{chan}}$.\\ 
In the ODE algorithm \cite{rid}, we need to adapt the time step $h$ when the parameter $N_{\text{chan}}$ varies, otherwise, we do not simulate the expected trajectories. Thinning algorithm in the \textit{channel model} speeds up the simulation by a factor 3 compared to the ODE method whereas in the \textit{subunit model} the factor is approximately 1. Such a difference is explained by the fact that the ratio of the simulation times between the flows of the \textit{subunit} and the \textit{channel} (for thinning algorithm) is bigger than the ratio of the simulation times between the vector fields of the \textit{subunit} and the \textit{channel} (for ODE algorithm).\\
Despite the complexity of the \textit{optimal bound} compared to the two others, it is the most efficient
one to simulate both the \textit{channel model} and the \textit{subunit model}.

\subsubsection{Spiking times}

Bio-scientists believe that the timing of action potentials is one of the characteristics of the nervous system which carries the most of informations. It has been shown experimentally \cite{verven} that if a neuron is repeatedly stimulated by identical pulses, both the amplitude and the timing of the action potentials is variable. In the sequel we numerically compare the mean value of the spiking time of the \textit{subunit} and \textit{channel model} to the one of the deterministic Hodgkin-Huxley model.\\    
Let $(x_{t})$ be the \textit{subunit model} or the \textit{channel model} defined on a filtered probability space $(\Omega,\mathcal{F}, \mathcal{F}_{t}, \mathbb{P}_{x})$. We consider that the stimulation is a monophasic current which produces only one action potential within a given time window $[0,T]$ as in Figure \ref{fig:1}. We suppose that a spike occurs when the membrane potential exceeds a certain value noted $\nu$. Let $\mathcal{T}$ be the spiking time defined by
\[
\mathcal{T}=\inf\{t\in{[0,T]} : V_{t}\geq{\nu}\}
\]
We are interested in the numerical computation of the mean and the standard deviation of $\mathcal{T}$ as a function of the number of channels. 
For low values of the parameters $N_{\text{Na}}$ and $N_{\text{K}}$ a spike may never occur. In this case, $\mathcal{T}=T$ and we do not count these trajectories in the Monte Carlo procedure. Thus, we evaluate the mean value of the spiking time conditionally on having a spike, $\mathbb{E}[\mathcal{T}|\mathcal{T}<T]$, with the following estimator $I_{M}=(1/M)\sum_{k=1}^{M}\mathcal{T}_{k}$ where $(\mathcal{T}_{k})$ are iid realizations of $\mathcal{T}$ conditionally on $\{\mathcal{T}<T\}$.

\begin{figure}[h!]
\includegraphics[width=.75\textwidth]{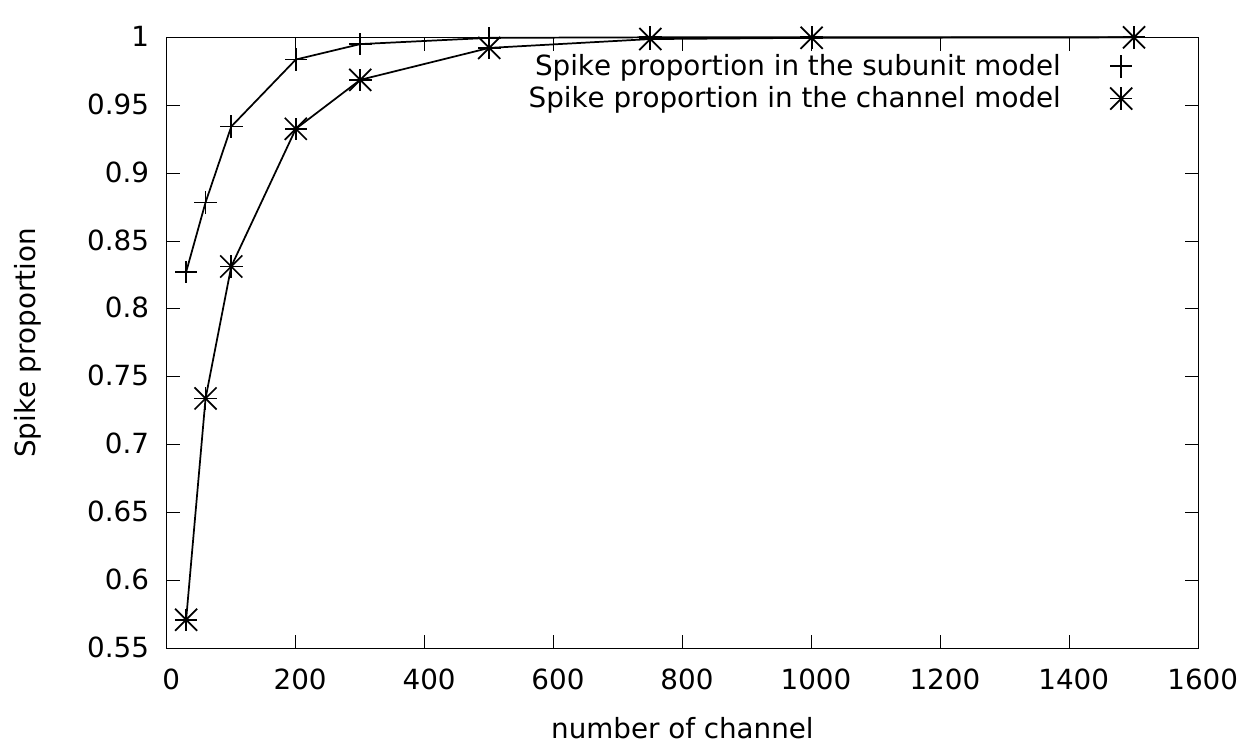}
\caption{The graph represents the proportion of spikes obtained with the \textit{subunit model} and the \textit{channel model} as a function of the number of channels $N_{\text{chan}}$.}

\label{fig:5}
\end{figure}
\begin{figure}[h!]
  \begin{subfigure}[b]{.5\linewidth}
    \includegraphics[width=\textwidth]{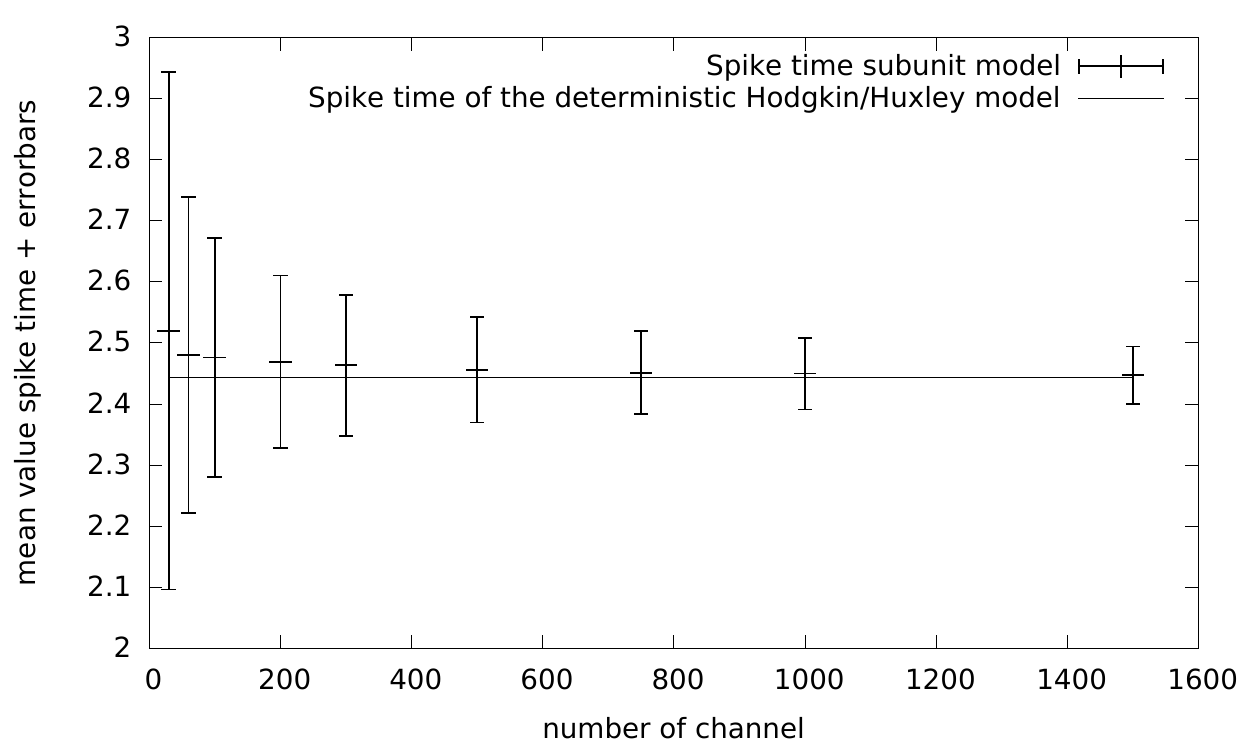}
  \end{subfigure} 
  \begin{subfigure}[b]{.5\linewidth}
    \includegraphics[width=\textwidth]{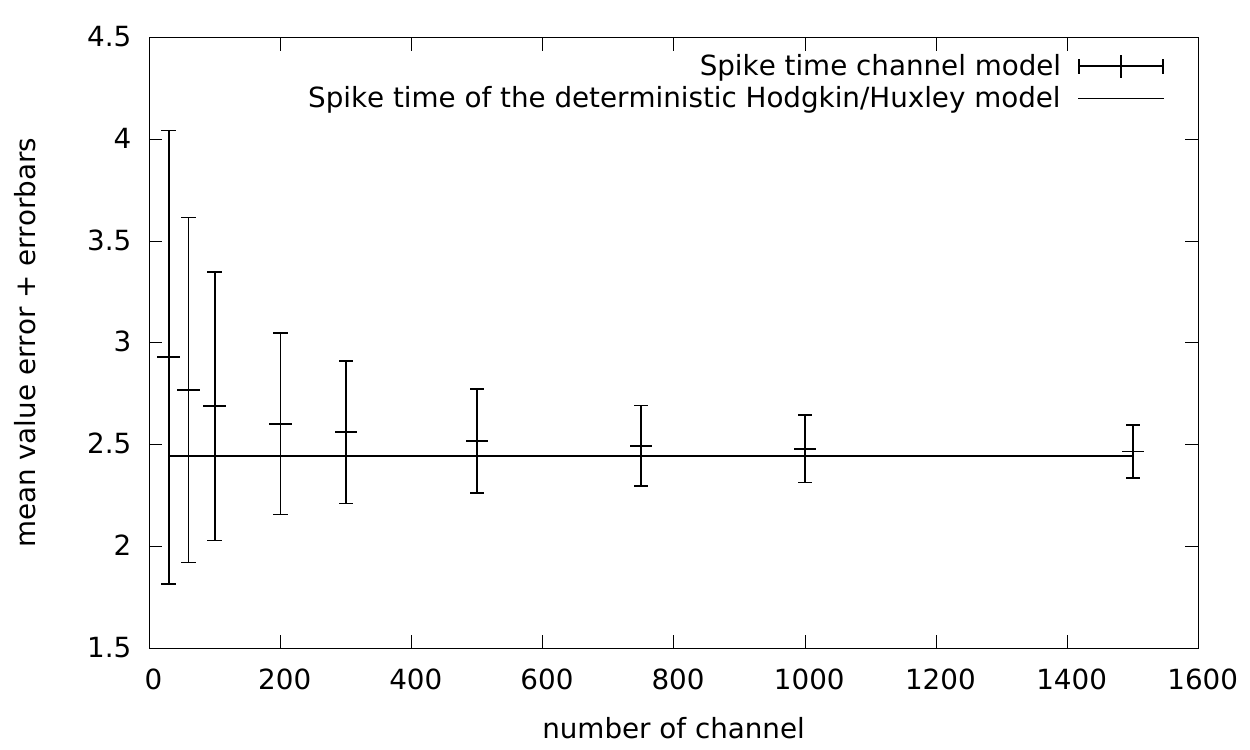}
  \end{subfigure}
  \caption{Mean value of the spiking time (ms) with standard deviation as a function of the number of channels $N_{\text{chan}}$. Left: \textit{subunit model}. Right: \textit{channel model}.}
\label{fig:6}
\end{figure}

It has been shown in \cite{wainfluid} that the deterministic limits of both the \textit{subunit} (Hodgkin-Huxley of dimension four \cite{hodgkin}) and the \textit{channel model} (Hodgkin-Huxley of dimension fourteen \cite{wainfluid}) are equivalent when the initial conditions satisfy a combinatorial relationship. We consider that, at time $t=0$, all the gates of the \textit{subunit model} are closed and all the channels of the \textit{channel model} are in the corresponding state, i.e state $\{m_{0}h_{0}\}$ for the sodium and $\{n_{0}\}$ for the potassium. These initial conditions satisfy the combinatorial relationship in \cite{wainfluid}. The initial conditions of both deterministic Hodgkin-Huxley models are also chosen so that they satisfy the binomial relation. Thus, the spiking time of these deterministic models is the same.
In the simulations, we take $T=10$, $\nu=60$, we consider that the stimulation is given by $I(t)=30\textbf{1}_{[1,2]}(t)$ and that $N_{\text{Na}}=N_{\text{K}}=N_{\text{chan}}$. In this case, the spiking time of the deterministic model is $\mathcal{T}^{\text{deter}}=2,443$.\\
Figure 6 illustrates the convergence of the mean spiking time of both the \textit{subunit} and the \textit{channel model} when the number of channels goes to infinity. For $N_{\text{chan}}=1500$ we see that the dispersion of the spiking time around its deterministic limit is approximately of order $10^{-1}$ ms for the \textit{subunit model} and of order $10^{-2}$ ms for the \textit{channel model}. Thus, a membrane patch with a number of channels superior to 1500 mimics the behaviour of the deterministic Hodgkin-Huxley model. For a number of channels inferior to 500, we see from figure 5 that the neuron may not respond to the stimuli. In this case, the dispersion of the spiking time ranges from approximately $10^{-1}$ and almost $1$ ms and are consistent with the observations in \cite{verven}. Since the simulation is exact the estimator $I_{M}$ is unbiased and errors due to the Monte Carlo procedure are of order of $M^{-1/2}$.

\begin{appendix}
\section{Rate functions and parameters}\label{parameters}
\[
\begin{array}{lll}
\alpha_{n}(V)=\frac{(0.1-0.01V)}{\exp(1-0.1V)-1}, & \alpha_{m}(V)=\frac{(2.5-0.1V)}{\exp(2.5-0.1V)-1}, &  \alpha_{h}(V)=0.07\exp(-\frac{V}{20}), \\
\vspace{0.1cm}\\
\beta_{n}(V)=0.125\exp(-\frac{V}{80}), & \beta_{m}(V)=4\exp(-\frac{V}{18}), & \beta_{h}(V)=\frac{1}{\exp(3-0.1V)+1}, \\
\end{array}\]
\[\begin{array}{lllllll}
V_{\text{Na}}=115, & g_{\text{Na}}=120, &
V_{\text{K}}=-12, & g_{\text{K}}=36, &
V_{\text{L}}=0, & g_{\text{L}}=0.3, &
C=1.
\end{array}
\]

\section{Markov schemes for the channel model}\label{Na_K}
Sodium (Na) scheme:\\
\[
\begin{array}{lllllll}
\textbf{m}_{\textbf{0}}\textbf{h}_{\textbf{0}} & \begin{array}{l} 3\alpha_{m}\\\longrightarrow\\\longleftarrow \\ \beta_{m} \end{array} & \textbf{m}_{\textbf{1}}\textbf{h}_{\textbf{0}} & \begin{array}{l} 2\alpha_{m}\\\longrightarrow\\\longleftarrow\\2\beta_{m} \end{array} & \textbf{m}_{\textbf{2}}\textbf{h}_{\textbf{0}} & \begin{array}{l} \alpha_{m}\\\longrightarrow\\\longleftarrow\\3\beta_{m} \end{array} & \textbf{m}_{\textbf{3}}\textbf{h}_{\textbf{0}}\\
\begin{array}{lll}\beta_{h}\uparrow\downarrow\alpha_{h} \end{array} & & \begin{array}{lll}\beta_{h}\uparrow\downarrow\alpha_{h} \end{array}  & & \begin{array}{lll}\beta_{h}\uparrow\downarrow\alpha_{h} \end{array}& & \begin{array}{lll}\beta_{h}\uparrow\downarrow\alpha_{h} \end{array}\\
 \textbf{m}_{\textbf{0}}\textbf{h}_{\textbf{1}} & \begin{array}{l} 3\alpha_{m}\\\longrightarrow\\\longleftarrow \\ \beta_{m} \end{array} & \textbf{m}_{\textbf{1}}\textbf{h}_{\textbf{1}} &\begin{array}{l} 2\alpha_{m}\\\longrightarrow\\\longleftarrow\\2\beta_{m} \end{array} & \textbf{m}_{\textbf{2}}\textbf{h}_{\textbf{1}} &\begin{array}{l} \alpha_{m}\\\longrightarrow\\\longleftarrow\\3\beta_{m} \end{array} & \textbf{m}_{\textbf{3}}\textbf{h}_{\textbf{1}}  \\
\end{array}
\]
\vspace{1cm}\\
Potassium (K) scheme:
\[
\begin{array}{lllllllll}
\textbf{n}_\textbf{0} & \begin{array}{l} 4\alpha_{n}\\\longrightarrow\\\longleftarrow \\ \beta_{n} \end{array} & \textbf{n}_\textbf{1} & \begin{array}{l} 3\alpha_{n}\\\longrightarrow\\\longleftarrow \\ 2\beta_{n} \end{array} & \textbf{n}_\textbf{2} & \begin{array}{l} 2\alpha_{n}\\\longrightarrow\\\longleftarrow \\ 2\beta_{n} \end{array} &\textbf{n}_\textbf{3} & \begin{array}{l} \alpha_{n}\\\longrightarrow\\\longleftarrow \\ 4\beta_{n} \end{array} & \textbf{n}_\textbf{4}
\end{array}
\]

\section{Rate of acceptance for the thinning of Poisson processes}

Let $N$ and $\tilde{N}$ be two Poisson processes with jump rate $\lambda$ and $\tilde{\lambda}$ respectively and jump times $(T_{n})_{n\geq1}$ and $(\tilde{T}_{n})_{n\geq1}$ respectively. Assume that $N$ is the thinning of $\tilde{N}$. Since $\mathbb{P}(\tilde{N}_{t}=0)=e^{-\int_{0}^{t}\tilde{\lambda}(s)ds}$, we define the rate of acceptance by $\mathbb{E}[N_{t}/\tilde{N}_{t}|\tilde{N}_{t}\geq1]$. In the case of Poisson processes this indicator takes the following form
\[\mathbb{E}[\frac{N_{t}}{\tilde{N}_{t}}|\tilde{N}_{t}\geq1]=\frac{\int_{0}^{t}\lambda(s)ds}{\int_{0}^{t}\tilde{\lambda}(s)ds}\tag{7}\]
To get $(7)$, we use the following result which is similar to the $n$-uplet of non-ordering uniform variables in the Poisson homogeneous case
\[f_{(\tilde{T}_{1},\ldots,\tilde{T}_{n}|\tilde{N}_{t}=n)}(t_{1},\ldots,t_{n})=\frac{\tilde{\lambda}(t_{1})\ldots
\tilde{\lambda}(t_{n})}{\left(\int_{0}^{t}\tilde{\lambda}(s)ds\right)^n}
\textbf{1}_{(t_{1},\ldots,t_{n})\in{[0,t]^{n}}}\tag{8}\] 
Equation $(8)$ gives an explicit formula of the conditional density of the vector $(\tilde{T}_{1},\ldots,\tilde{T}_{n}|\tilde{N}_{t}=n)$. Note that we do not consider any ordering in points $(\tilde{T}_{k})_{0\leq{k}\leq{n}}$ and that conditionally to $\tilde{N}_{t}=n$, the points $\tilde{T}_{1},\ldots,\tilde{T}_{n}$ are independent with density $\left(\tilde{\lambda}(s)/\int_{0}^{t}\tilde{\lambda}(u)du\right)\textbf{1}_{s\in{[0,t]}}$.\\
With $(8)$ one is able to determine that 
\[\mathcal{L}(N_{t}|\tilde{N}_{t}=n)=\mathcal{B}(n,p)\tag{9}\]
with $p=\int_{0}^{t}\lambda(s)ds/\int_{0}^{t}\tilde{\lambda}(s)ds$ by noting that for $k\leq{n}$ \[\{N_{t}=k|\tilde{N}_{t}=n\}=\bigcup_{1\leq{i_{1}}<\ldots<i_{k}\leq{n}}\Big{[}
\bigcap_{i\in{\{i_{1},\ldots,i_{k}\}}}\{U_{i}\leq{\frac{\lambda}{\tilde{\lambda}}(\tilde{T}_{i})|\tilde{N}_{t}=n}\}
\bigcap_{i\in{\{i_{1},\ldots,i_{k}\}}^{c}}\{U_{i}>\frac{\lambda}{\tilde{\lambda}}(\tilde{T}_{i})|\tilde{N}_{t}=n\}\Big{]}\]
and then
\[
\mathbb{P}(N_{t}=k|\tilde{N}_{t}=n)=\left(\begin{array}{c}n\\k\end{array}\right)\mathbb{P}\left(U_{i}\leq{\frac{\lambda}{\tilde{\lambda}}(\tilde{T}_{i})}|\tilde{N}_{t}=n\right)^{k}
\mathbb{P}\left(U_{i}>\frac{\lambda}{\tilde{\lambda}}(\tilde{T}_{i})|\tilde{N}_{t}=n\right)^{n-k}
\]
where $(U_{i})$ are independent variables uniformly distributed in $[0,1]$, independent of $(\tilde{T}_{i})$.
Thus, the law of the number of accepted points is binomial conditionally on the number of proposed points. Then, we find (7) by noting that
\[
\mathbb{E}[\frac{N_{t}}{\tilde{N}_{t}}|\tilde{N}_{t}\geq1]=\frac{1}{\mathbb{P}(\tilde{N}_{t}\geq1)}\sum_{n\geq1}\frac{1}{n}\mathbb{E}[N_{t}|\tilde{N}_{t}=n]\mathbb{P}(\tilde{N}_{t}=n)
\]
We should note that to find the explicit formula (3), we work almost exclusively on the process $(\tilde{N}_{t})$. The only characteristic of $N$ that we use is its jump rate $\lambda$ and we use it to evaluate the thinning probabilities $\left(\lambda/\tilde{\lambda}\right)(.)$.

\end{appendix}

\bibliographystyle{plain}
\bibliography{bibPDMP}

\begin{thebibliography}{10}

\bibitem{alfonsi}
A.~Alfonsi, E.~Cancés, G.~Turinci, B.~Di~Ventura, and W.~Huisinga.
\newblock Adaptive hybrid simulation of hybrid stochastic and deterministic
  models for biochemical reactions.
\newblock {\em ESAIM: proceedings}, 14:1--13, 2005.

\bibitem{ander}
D.F. Anderson, B.~Ermentrout, and P.J. Thomas.
\newblock Stochastic representation of ion channel kinetics and exact
  stochastic simulation of neuronal dynamics.
\newblock {\em Journal of {C}omputational {N}euroscience}, 38:67--82, 2015.

\bibitem{bouguet}
F.~Bouguet.
\newblock Quantitative speeds of convergence for exposure to food contaminants.
\newblock {\em ESAIM: Probability and Statistics}, 19:482--501, 2015.

\bibitem{stried}
E.~Buckwar and M.G. Riedler.
\newblock An exact stochastic hybrid model of excitable membranes including
  spatio-temporal evolution.
\newblock {\em Journal of Mathematical Biology}, 63:1051--1093, 2011.

\bibitem{malri}
D.~Chafaï, F.~Malrieu, and K.~Paroux.
\newblock On the long time behavior of the tcp window size process.
\newblock {\em Stochastic Processes and their Applications}, 120:1518–1534,
  2010.

\bibitem{chowwhite}
C.C. Chow and J.A. White.
\newblock Spontaneous action potentials due to channel fluctuations.
\newblock {\em Biophysical Journal}, 71:3013--3021, 1996.

\bibitem{defelice}
J.R. Clay and De{F}elice.
\newblock Relationship between membrane excitability and single channel
  open-close kinetics.
\newblock {\em Biophysical Journal}, 42:151--157, 1983.

\bibitem{cocofiab}
C.~Cocozza-Thivent.
\newblock {\em Processus stochastiques et fiabilité des systèmes}.
\newblock Springer, 1997.

\bibitem{davarticle}
M.H.A. Davis.
\newblock Piecewise-deterministic {M}arkov processes: A general class of
  non-diffusion stochastic models.
\newblock {\em Journal of the Royal statistical Society}, 46:353--388, 1984.

\bibitem{dav}
M.H.A. Davis.
\newblock {\em Markov Models and Optimization}.
\newblock Chapman and Hall, {L}ondon, 1993.

\bibitem{dev}
L.~Devroye.
\newblock {\em Non-uniform random variate generation}.
\newblock Springer-Verlag, New York Inc., 1986.

\bibitem{doumic}
M.~Doumic, M.~Hoffmann, N.~Krell, and L.~Robert.
\newblock Statistical estimation of a growth-fragmentation model observed on a
  genealogical tree.
\newblock {\em Bernoulli Society for Mathematical Statistics and Probability},
  21:1760--1799, 2015.

\bibitem{fox}
R.F. Fox.
\newblock Stochastic {V}ersions of the {H}odgkin-{H}uxley {E}quations.
\newblock {\em Biophysical {J}ournal}, 72:2068--2074, 1997.

\bibitem{gold}
J.H. Goldwyn, N.S. Imennov, M.~Famular, and E.~Shea-Brown.
\newblock Stochastic differential equation models for ion channels noise in
  {H}odgkin-{H}uxley neurons.
\newblock {\em Physical Review E: Statistical, Nonlinear, and Soft Matter
  Physics}, 83 (4 {P}t 1), 2011.

\bibitem{Bertil}
B.~Hille.
\newblock {\em Ionic channels of excitable membranes}.
\newblock Sinauer Associates Inc, 1992.

\bibitem{hodgkin}
A.L. Hodgkin and A.F. Huxley.
\newblock A quantitative description of membrane current and its application to
  conduction and excitation in nerve.
\newblock {\em Journal of Physiology}, 117:500--544, 1952.

\bibitem{kalrand}
O.~Kallenberg.
\newblock {\em Random measures}.
\newblock Akademie-Verlag Berlin, 1975.

\bibitem{kal}
O.~Kallenberg.
\newblock {\em Foundations of Modern Probability}.
\newblock Springer,New York, 2002.

\bibitem{lew}
P.A.W. Lewis and G.S. Shedler.
\newblock Simulation of nonhomogeneous {P}oisson processes by thinning.
\newblock {\em Naval Research Logistics Quarterly}, 26:403--413, 1979.

\bibitem{mino}
H.~Mino, J.T. Rubinstein, and J.A White.
\newblock Comparison of algorithms for the simulation of action potentials with
  stochastic sodium channels.
\newblock {\em Annals of Biomedical Engineering}, 30:578--587, 2002.

\bibitem{lecar}
C.~Morris and H.~Lecar.
\newblock Voltage oscillations in the barnacle giant muscle fiber.
\newblock {\em Biophysical Journal}, 35:193--213, 1981.

\bibitem{ogata}
Y.~Ogata.
\newblock {O}n {L}ewis’ {S}imulation {M}ethod for {P}oint {P}rocesses.
\newblock {\em {I}{E}{E}{E} Transactions on Information Theory}, 27, 1981.

\bibitem{soudry}
P.~Orio and D.~Soudry.
\newblock Simple, fast and accurate implementation of the diffusion
  approximation algorithm for stochastic ion channels with multiple states.
\newblock {\em PLoS one}, 7, 2012.

\bibitem{wainfluid}
K.~Pakdaman, M.~Thieullen, and G.~Wainrib.
\newblock Fluid limit theorems for stochastic hybrid systems with application
  to neuron models.
\newblock {\em advanced applied probability}, 42:761--794, 2010.

\bibitem{rid}
M.G. Riedler.
\newblock Almost sure convergence of numerical approximations for {P}iecewise
  {D}eterministic {M}arkov {P}rocesses.
\newblock {\em Journal of Computational and Applied Mathematics}, pages 50--71,
  2012.

\bibitem{rubin}
J.T. Rubinstein.
\newblock Threshold fluctuations in an {N} sodium channel model of the node of
  ranvier.
\newblock {\em Biophysical Journal}, 68:779--785, 1995.

\bibitem{skaugen}
E.~Skaugen and L.~Walloe.
\newblock Firing behaviour in a stochastic nerve membrane model based upon the
  {H}odgkin-{H}uxley equations.
\newblock {\em Acta Physiologica}, 107:343--363, 1979.

\bibitem{veltz}
R.~Veltz.
\newblock A new twist for the simulation of hybrid systems using the true jump
  method.
\newblock {\em arXiv[math]}, 2015.

\bibitem{verven}
A.A. Verveen and H.E. Derksen.
\newblock Fluctuation phenomena in nerve mebrane.
\newblock {\em proceedings of the IEEE}, 56:906--916, 1968.

\end{thebibliography}

\end{document}